\documentclass[aap,MSNbibl,dvips]{arximspdf}
\usepackage{mathbh}

\doi{10.1214/13-AAP983} 
\volume{24}
\issue{6}
\pubyear{2014}
\firstpage{2491}
\lastpage{2526}

\makeatletter
\newcommand{\widebar}{\overline}

\newcommand{\MD}{\mathrm{MD}}
\newcommand{\WL}{\mathrm{WL}}
\newcommand{\rrvert}{\vert}
\newcommand{\llvert}{\vert}
\newtheorem{theorem}{Theorem}[section]
\newtheorem{lemma}[theorem]{Lemma}
\newtheorem{proposition}[theorem]{Proposition}
\newproclaim{remark}{Remark}
\makeatother

\begin{document}
\begin{frontmatter}

\title{Limit theorems for nondegenerate \textit{U}-statistics of continuous
semimartingales}
\runtitle{$U$-statistics of semimartingales}

\begin{aug}
\author[a]{\fnms{Mark} \snm{Podolskij}\ead[label=e1]{m.podolskij@uni-heidelberg.de}},
\author[a]{\fnms{Christian} \snm{Schmidt}\ead[label=e2]{christian-schmidt@uni-heidelberg.de}}
\and
\author[b]{\fnms{Johanna F.} \snm{Ziegel}\corref{}\ead[label=e3]{johanna.ziegel@stat.unibe.ch}}
\runauthor{M. Podolskij, C. Schmidt and J. F. Ziegel}
\affiliation{Heidelberg University, Heidelberg University and
University of Bern}
\address[a]{M. Podolskij\\
C. Schmidt\\
Department of Mathematics\\
Heidelberg University\\
INF 294\\
69120 Heidelberg\\
Germany\\
\printead{e1}\\
\phantom{E-mail:\ }\printead*{e2}} 
\address[b]{J. F. Ziegel\\
Department of Mathematics and Statistics\\
Institute of Mathematical Statistics\\
\quad and Actuarial Science\\
University of Bern\\
Sidlerstrasse 5\\
3012 Bern\\
Switzerland\\
\printead{e3}}
\end{aug}

\received{\smonth{11} \syear{2012}}
\revised{\smonth{10} \syear{2013}}

%
\begin{abstract}
This paper presents the asymptotic theory for nonde\-gene\-rate
$U$-statistics of high
frequency observations of continuous It\^o semimartingales. We prove uniform
convergence in probability and show a functional stable central limit
theorem for the standardized version
of the $U$-statistic. The limiting process in the central limit theorem
turns out to be conditionally Gaussian with
mean zero. Finally, we indicate potential statistical applications of
our probabilistic results.
\end{abstract}

%
\begin{keyword}[class=AMS]
\kwd[Primary ]{60F05}
\kwd{60F15}
\kwd{60F17}
\kwd[; secondary ]{60G48}
\kwd{60H05}
\end{keyword}
\begin{keyword}
\kwd{High frequency data}
\kwd{limit theorems}
\kwd{semimartingales}
\kwd{stable convergence}
\kwd{$U$-statistics}
\end{keyword}

\end{frontmatter}

\section{Introduction} \label{Intro}

Since the seminal work by Hoeffding \cite{Hoe}, $U$-statistics have been
widely investigated by probabilists
and statisticians. Nowadays, there exists a vast amount of literature
on the asymptotic properties of $U$-statistics in the case of independent
and identically distributed (i.i.d.) random variables or in the
framework of weak dependence. We refer to \cite
{KoroljukBorovskich1994} for a comprehensive account of the asymptotic theory
in the classical setting. In \cite{BBD1,BBD2,DK}, the authors treat
limit theorems for $U$-statistics under various mixing conditions,
while the corresponding theory for long memory processes has been
studied, for example, in \cite{DT,GT}; see \cite{HW}
for a recent review of the properties of $U$-statistics in various
settings. The most powerful tools for proving asymptotic
results for $U$-statistics include the classical Hoeffding decomposition
(see, e.g.,~\cite{Hoe}), Hermite expansions (see, e.g.,~\mbox{\cite{DT,DTb}})
and the empirical process approach; see, for example,~\cite{BZ}.
Despite the activity of this field of research, $U$-statistics for high
frequency observations of a time-continuous process
have not been studied in the literature thus far. The notion of
high frequency data refers to the sampling scheme in which the time
step between two consecutive observations converges
to zero while the time span remains fixed. This concept is also known
under the name of infill asymptotics. Motivated by the prominent role
of semimartingales in mathematical finance, in this paper we present novel
asymptotic results for high frequency observations of It\^o
semimartingales and demonstrate some statistical applications.

The seminal work of Jacod \cite{J3} marks the starting point for
stable limit theorems for semimartingales.
Stimulated by the increasing popularity of semimartingales as natural
models for asset pricing, the asymptotic theory
for partial sums processes of continuous and discontinuous It\^o semimartingales
has been developed in \cite{BGJPS,J2,KP}; see also the recent book
\cite{JP}. We refer to \cite{PV} for a short survey of limit theorems
for semimartingales. More recently, asymptotic theory
for It\^o semimartingales observed with errors has been investigated
in \cite{JPV}.

The methodology we employ to derive a limit theory for $U$-statistics of
continuous It\^o semimartingales is an intricate combination and
extension of some of the techniques developed in the series of papers
mentioned in the previous paragraph and the empirical process approach
to $U$-statistics.

In this paper we consider a one-dimensional continuous It\^o
semimartingale of the form
\[
X_t= x + \int_0^t a_s
\,ds + \int_0^t \sigma_s
\,dW_s, \qquad t\geq0,
\]
defined on a filtered probability space $(\Omega,\mathcal
{F},(\mathcal F_t)_{t\geq0},\mathbb{P})$ (which satisfies the usual
assumptions), where $x \in\mathbb{R}$, $(a_s)_{s\geq0}$,
$(\sigma_s)_{s\geq0}$ are stochastic processes, and $W$ is a standard
Brownian motion. The underlying observations of $X$
are
\[
X_{{i}/{n}}, \qquad i=0, \ldots, [nt],
\]
and we are in the framework of infill asymptotics, that
is,~$n\rightarrow\infty$. In order to present our main results,
we introduce some notation. We define
\begin{eqnarray*}
\mathcal{A}_t^n(d) &:=& \bigl\{\mathbf{i}=(i_1,
\ldots,i_d) \in\mathbb{N}^d\dvtx 1\leq i_1<i_2<
\cdots<i_d\leq[nt]\bigr\},
\\
Z_\mathbf{s}&:=&(Z_{s_1}, \ldots, Z_{s_d}),
\qquad\mathbf{s}\in\mathbb{R}^d,
\end{eqnarray*}
where $Z=(Z_t)_{t \in\mathbb{R}}$ is an arbitrary stochastic process.
For any continuous function $H\dvtx\mathbb{R}^d \rightarrow\mathbb{R}$,
we define the $U$-statistic $U(H)_t^n$ of order $d$
as
%
%
\begin{equation}
\label{ustatistic} U(H)_t^n = {\pmatrix{n
\cr
d}}^{-1} \sum_{\mathbf{i}\in\mathcal{A}_t^n(d)} H\bigl( \sqrt{n}
\Delta_{\mathbf{i}}^n X\bigr)
\end{equation}
with $\Delta_{\mathbf{i}}^n X = X_{\mathbf{i}/n} - X_{(\mathbf{i}-1)/n}$.
For a multi-index $\mathbf{i} \in\mathbb{N}^d$, the vector $\mathbf
{i} - 1$ denotes the multi-index obtained by componentwise subtraction
of $1$ from $\mathbf{i}$.
In the following we assume that the function $H$ is symmetric, that
is, for all $x=(x_1,\ldots,x_d)\in\mathbb{R}^d$ and all permutations
$\pi$ of
$\{1, \ldots,d\}$, it holds that $H(\pi x)=H(x)$, where $\pi x =
(x_{\pi(1)},\ldots,x_{\pi(d)})$.

Our first result determines the asymptotic behavior of $U(H)_t^n$,
\[
U(H)_t^n \stackrel{\mathrm{u.c.p.}} {
\longrightarrow}U(H)_t:= \int_{[0,t]^d}
\rho_{\sigma_{\mathbf{s}}} (H) \,d\mathbf{s},
\]
where $Z^n \stackrel{\mathrm{u.c.p.}}{\longrightarrow}Z$ denotes
uniform convergence in probability, that is,
for any \mbox{$T>0$}, $\sup_{t\in[0,T]} |Z_t^n - Z_t|$ $\stackrel{\mathbb
{P}}{\longrightarrow}0$, and
%
%
\begin{equation}
\label{rho} \rho_{\sigma_{\mathbf{s}}} (H):= \int_{\mathbb{R}^d} H(
\sigma_{s_1} u_1, \ldots, \sigma_{s_d}
u_d) \varphi_d(\mathbf{u}) \,d\mathbf{u}
\end{equation}
with $\varphi_d$ denoting the density of the $d$-dimensional standard
Gaussian law $\mathcal N_d(0,\mathbf{I}_d)$. The second result of this
paper is the stable functional central limit
theorem
\[
\sqrt{n} \bigl(U(H)^n - U(H)\bigr) \stackrel{\mathrm{st}} {\longrightarrow}L,
\]
where $\stackrel{\mathrm{st}}{\longrightarrow}$ denotes stable convergence in
law, and the function $H$ is assumed to be even
in each coordinate. The limiting process $L$ lives on an extension of
the original
probability space $(\Omega,\mathcal{F},(\mathcal F_t)_{t\geq
0},\mathbb{P})$ and it turns out to be Gaussian with mean zero
conditionally on the original $\sigma$-algebra
$\mathcal F$. The proofs of the asymptotic results rely upon a
combination of recent limit theorems for semimartingales
(see, e.g.,~\cite{J3,JP,KP}) and empirical processes techniques.

The paper is organized as follows. In Section~\ref{large-numbers} we
present the law of large numbers for the $U$-statistic $U(H)_t^n$. The
associated functional stable central limit \mbox{theorem} is provided in
Section~\ref{secSCLT}. Furthermore, we derive a standard central
limit theorem in Section~\ref{seccondvar}. In Section~\ref{secstatex}
we demonstrate statistical applications of our limit theory including
Gini's mean difference, homoscedasticity testing
and Wilcoxon statistics for testing of structural breaks. Some
technical parts of the proofs are deferred to Section~\ref{proofs}.

\section{Preliminaries}\label{secprelim}
We consider the continuous diffusion model
%
%
\begin{equation}
\label{diffusion} X_t= x + \int_0^t
a_s \,ds + \int_0^t
\sigma_s \,dW_s, \qquad t\geq0,
\end{equation}
where\vspace*{1pt} $(a_s)_{s\geq0}$ is a c\`agl\`ad process, $(\sigma_s)_{s\geq
0}$ is a
c\`adl\`ag process, both adapted to the filtration $(\mathcal F_s)_{s
\ge0}$. Define the functional class $C_p^k(\mathbb{R}^d)$ via
\begin{eqnarray*}
C_p^k\bigl(\mathbb{R}^d \bigr)&:=&\bigl\{f\dvtx
\mathbb{R}^d \rightarrow\mathbb{R} |  f\in C^k\bigl(
\mathbb{R}^d \bigr) \mbox{ and all derivatives up to order $k$}
\\
&&\hspace*{143pt} \mbox{are of polynomial growth}\bigr\}.
\end{eqnarray*}
Note\vspace*{-1pt} that $H\in C^0_p(\mathbb{R}^d )$ implies that $\rho_{\sigma
_{\mathbf{s}}} (H)<\infty$ almost surely. For any vector
$y\in\mathbb{R}^d$, we denote by $\|y\|$ its maximum norm; for any
function $f\dvtx\mathbb{R}^d \rightarrow\mathbb{R}$, $\|f\|_{\infty}$
denotes its supremum norm. Finally, for any $z\neq0$, $\Phi_z$ and
$\varphi_z$ stand for the distribution function and density of the
Gaussian law $\mathcal N(0,z^2)$, respectively; $\Phi_0$ denotes the
Dirac measure at the origin.
The bracket $[M,N]$ denotes the covariation process of two local martingales
$M$ and $N$.

\section{Law of large numbers}\label{large-numbers}
We start with the law of large numbers, which describes the limit of
the $U$-statistic $U(H)_t^n$ defined at (\ref{ustatistic}).
First of all, we remark that the processes $(a_s)_{s\geq0}$ and
$(\sigma_{s-})_{s\geq0}$ are locally bounded, because they
are both c\`agl\`ad. Since the main results of this subsection
(Proposition~\ref{prop1} and Theorem~\ref{th1}) are
\textit{stable under stopping}, we may assume without loss of
generality that
%
%
\begin{equation}
\label{probound1} \mbox{The processes $a$ and $\sigma$ are bounded in
$(\omega,t)$.}
\end{equation}
A detailed justification of this statement can be found in \cite
{BGJPS}, Section~3.

We start with the representation of the process $U(H)_t^n$ as an
integral with respect to a certain empirical
random measure. For this purpose let us introduce the quantity
%
%
\begin{equation}
\label{alpha} \alpha_j^n:= \sqrt{n}
\sigma_{({j-1})/{n}} \Delta_j^n W, \qquad j\in
\mathbb{N},
\end{equation}
which serves as a first order approximation of the increments $\sqrt
{n} \Delta_j^n X$. The empirical distribution
function associated with the random variables $(\alpha
_j^n)_{1\leq j\leq[nt]}$ is defined as
%
%
\begin{equation}
\label{fn} F_n(t,x):= \frac{1}n \sum
_{j=1}^{[nt]} \mathbh{1}_{\{\alpha_j^n\leq
x\}}, \qquad x\in
\mathbb{R}, t\geq0.
\end{equation}
Notice that, for any fixed $t\geq0$, $F_n(t, \cdot)$ is a finite
random measure.
Let $\widetilde{U}(H)_t^n$ be the $U$-statistic based on $\alpha
_j^n$'s, that is,
%
%
\begin{equation}
\label{tildeustatistic} \widetilde{U}(H)_t^n = {\pmatrix{n
\cr
d}}^{-1} \sum_{\mathbf{i}\in
\mathcal{A}_t^n(d)} H\bigl(
\alpha_{\mathbf{i}}^n\bigr).
\end{equation}
The functional $U_t^{\prime n}(H)$ defined as
%
%
\begin{equation}
\label{uprime} U_t^{\prime n}(H):= \int_{\mathbb{R}^d}
H(\mathbf x) F_n^{\otimes d} (t, d\mathbf x),
\end{equation}
where
\[
F_n^{\otimes d} (t, d\mathbf x):= F_n(t,
dx_1)\cdots F_n(t, dx_d),
\]
is closely related to the process $\widetilde{U}(H)_t^n$; in fact, if
both are written out as multiple sums over nondecreasing multi-indices,
then their summands coincide on the set $\mathcal{A}_t^n(d)$. They
differ for multi-indices that have at least two equal components.
However, the number of these diagonal multi-indices is of order
$O(n^{d-1})$. We start with a simple lemma,
which we will often use throughout the paper. We omit a formal proof
since it follows by standard arguments.

%
\begin{lemma} \label{vague}
Let $Z_n,Z\dvtx[0,T] \times\mathbb{R}^m \rightarrow\mathbb{R}$, $n \ge
1$, be random positive functions such that $Z_n(t,\cdot)$ and
$Z(t,\cdot)$
are finite random measures on $\mathbb{R}^m$ for any $t\in[0,T]$.
Assume that
\[
Z_n(\cdot, \mathbf x) \stackrel{\mathrm{u.c.p.}} {\longrightarrow
}Z(\cdot, \mathbf x),
\]
for\vspace*{1pt} any fixed $\mathbf x\in\mathbb{R}^m$, and
$\sup_{t\in[0,T], \mathbf x\in\mathbb{R}^m} Z(t, \mathbf x)$, $\sup
_{t\in[0,T], \mathbf x\in\mathbb{R}^m} Z_n(t, \mathbf x)$, $n \ge
1$, are bounded random variables.
Then, for any continuous function $Q\dvtx \mathbb{R}^m \rightarrow
\mathbb
{R}$ with compact support,
we obtain that
\[
\int_{\mathbb{R}^m} Q(\mathbf x) Z_n(\cdot, d\mathbf x)
\stackrel{\mathrm{u.c.p.}} {\longrightarrow}\int_{\mathbb{R}^m} Q(
\mathbf x) Z(\cdot, d\mathbf x).
\]
\end{lemma}
The next proposition determines the asymptotic behavior of the
empirical distribution function
$F_n(t,x)$ defined at (\ref{fn}), and the $U$-statistic $U_t^{\prime
n}(H)$ given at (\ref{uprime}).

%
\begin{proposition} \label{prop1}
Assume that $H \in C_p^0(\mathbb{R}^d)$. Then, for any fixed $x\in
\mathbb{R}$,
it holds that
%
%
\begin{equation}
\label{fnconv} F_n(t,x) \stackrel{\mathrm{u.c.p.}} {
\longrightarrow}F(t,x):= \int_0^t
\Phi_{\sigma_s} (x) \,ds.
\end{equation}
Furthermore, we obtain that
%
%
\begin{equation}
\label{uprimeconv} U_t^{\prime n}(H) \stackrel{\mathrm{u.c.p.}} {
\longrightarrow}U(H)_t:= \int_{[0,t]^d}
\rho_{\sigma_{\mathbf
{s}}} (H) \,ds,
\end{equation}
where the quantity $\rho_{\sigma_{\mathbf{s}}} (H)$ is defined at
(\ref{rho}).
\end{proposition}
\begin{pf} Recall that we always assume (\ref{probound1}) without loss
of generality.
Here and throughout the paper, we denote by $C$ a generic positive
constant, which may change from
line to line; furthermore, we write $C_p$ if we want to emphasize the
dependence of $C$ on an external parameter $p$.
We first show the convergence in (\ref{fnconv}). Set $\xi_j^n:=n^{-1}
\mathbh{1}_{\{\alpha_j^n\leq x\}}$. It obviously
holds that
\[
\sum_{j=1}^{[nt]} \mathbb{E}\bigl[
\xi_j^n| \mathcal F_{(j-1)/{n}}\bigr] =
\frac{1}n \sum_{j=1}^{[nt]}
\Phi_{\sigma_{(j-1)/{n}}} (x) \stackrel{\mathrm{u.c.p.}}
{\longrightarrow} F(t,x),
\]
for any fixed $x\in\mathbb{R}$, due to Riemann integrability of the
process $\Phi_\sigma$. On the other hand,
we have for any fixed $x\in\mathbb{R}$,
\[
\sum_{j=1}^{[nt]} \mathbb{E}\bigl[\bigl|
\xi_j^n\bigr|^2| \mathcal F_{(j-1)/{n}}\bigr]
= \frac{1}{n^2} \sum_{j=1}^{[nt]}
\Phi_{\sigma_{(j-1)/{n}}} (x) \stackrel{\mathbb{P}} {\longrightarrow}0.
\]
This immediately implies the convergence (see \cite{JP}, Lemma 2.2.11, page~577)
\[
F_n(t,x) - \sum_{j=1}^{[nt]}
\mathbb{E}\bigl[\xi_j^n| \mathcal F_{(j-1)/{n}}\bigr]
= \sum_{j=1}^{[nt]} \bigl(
\xi_j^n- \mathbb{E}\bigl[\xi_j^n|
\mathcal F_{(j-1)/{n}}\bigr] \bigr) \stackrel{\mathrm{u.c.p.}} {
\longrightarrow}0,
\]
which completes the proof of (\ref{fnconv}). If $H$ is compactly
supported, then the convergence in (\ref{uprimeconv}) follows directly from
(\ref{fnconv}) and Lemma~\ref{vague}.

Now, let $H \in C_p^0(\mathbb{R}^d)$ be arbitrary. For any $k\in
\mathbb{N}$, let $H_k \in C_p^0(\mathbb{R}^d)$ be a function with
$H_k=H$ on
$[-k,k]^d$ and $H_k=0$ on $([-k-1,k+1]^d)^c$. We already know that
\[
U^{\prime n}(H_k) \stackrel{\mathrm{u.c.p.}} {
\longrightarrow}U(H_k),
\]
for any fixed $k$, and $U(H_k) \stackrel{\mathrm
{u.c.p.}}{\longrightarrow}U(H)$ as $k\rightarrow\infty$.
Since the function $H$ has polynomial growth,
that is, $|H(\mathbf x)| \leq C (1 + \|\mathbf x\|^q)$ for some $q>0$,
we obtain for any $p>0$
%
%
\begin{equation}
\label{pmoments} \mathbb{E}\bigl[\bigl|H\bigl(\alpha_{\mathbf{i}}^n
\bigr)\bigr|^p\bigr]\leq C_p \mathbb{E}\bigl[\bigl(1+ \bigl\|
\alpha_{\mathbf{i}}^n \bigr\|^{qp}\bigr)\bigr] \leq
C_p
\end{equation}
uniformly in $\mathbf{i}$, because the process $\sigma$ is bounded.
Statement (\ref{pmoments}) also holds
for~$H_k$. Recall that the function $H-H_k$ vanishes on $[-k,k]^d$. Hence,
we deduce by (\ref{pmoments}) and Cauchy--Schwarz inequality that
\begin{eqnarray*}
&& \mathbb{E}\Bigl[\sup_{t\in[0,T]}\bigl|U_t^{\prime n}(H-H_k)\bigr|
\Bigr]
\\
&&\qquad \leq C {\pmatrix{n
\cr
d}}^{-1} \sum_{1\le i_1, \ldots, i_d \leq[nT]}
\bigl(\mathbb{E}[\mathbh{1}_{\{|\alpha_{i_1}^n|\geq k\}} + \cdots+
\mathbh{1}_{\{|\alpha_{i_d}^n|\geq k \}}]
\bigr)^{1/2}
\\
&&\qquad \leq C_T \sup_{s\in[0,T]} \bigl(\mathbb{E}\bigl[1-
\Phi_{\sigma
_s}(k)\bigr] \bigr)^{1/2} \rightarrow0
\end{eqnarray*}
as $k\rightarrow\infty$. This completes the proof of (\ref{uprimeconv}).
\end{pf}

Proposition~\ref{prop1} implies the main result of this section.

%
\begin{theorem} \label{th1}
Assume that $H \in C_p^0(\mathbb{R}^d)$. Then it holds that
%
%
\begin{equation}
\label{LLN} U(H)_t^n \stackrel{\mathrm{u.c.p.}} {
\longrightarrow}U(H)_t:= \int_{[0,t]^d}
\rho_{\sigma_{\mathbf{s}}} (H) \,ds,
\end{equation}
where the quantity $\rho_{\sigma_{\mathbf{s}}} (H)$ is defined at
(\ref{rho}).
\end{theorem}
\begin{pf} In Section~\ref{proofs} we will show that
%
%
\begin{equation}
\label{sec3a} U(H)^n - \widetilde{U}(H)^n \stackrel{
\mathrm{u.c.p.}} {\longrightarrow}0,
\end{equation}
where the functional $\widetilde{U}(H)_t^n$ is given at (\ref
{tildeustatistic}). In view of Proposition~\ref{prop1}, it remains to
prove that $\widetilde{U}(H)_t^n - U_t^{\prime n}(H) \stackrel
{\mathrm{u.c.p.}}{\longrightarrow}0$. But due
to the symmetry of $H$ and estimation~(\ref{pmoments}),
we obviously obtain that
\[
\mathbb{E}\Bigl[\sup_{t\in[0,T]}\bigl|\widetilde{U}(H)_t^n
- U_t^{\prime
n}(H)\bigr|\Bigr] \leq\frac{C_T}{n} \rightarrow0,
\]
since the summands in $\widetilde{U}(H)_t^n$ and $U_t^{\prime n}(H)$
are equal except for diagonal multi-indices.
\end{pf}

%
\begin{remark} \label{rem1}
The result of Theorem~\ref{th1} can be extended to weighted $U$-statis\-tics of the type
%
%
\begin{equation}
\label{weightustatistic} U(H;X)_t^n:= {\pmatrix{n
\cr
d}}^{-1} \sum_{\mathbf{i}\in\mathcal
{A}_t^n(d)} H\bigl(X_{({\mathbf{i}-1})/{n}};
\sqrt{n} \Delta_{\mathbf{i}}^n X\bigr).
\end{equation}
Here, $H\dvtx\mathbb{R}^d \times\mathbb{R}^d \to\mathbb{R}$ is assumed
to be continuous and symmetric in the first and last $d$ arguments.
Indeed, similar methods of proof imply the u.c.p. convergence
\[
U(H;X)_t^n \stackrel{\mathrm{u.c.p.}} {
\longrightarrow}U(H;X)_t= \int_{[0,t]^d}
\rho_{\sigma_{\mathbf{s}}} (H;X_{\mathbf{s}}) \,d\mathbf{s},
\]
with
\[
\rho_{\sigma_{\mathbf{s}}} (H;X_{\mathbf{s}}):= \int_{\mathbb
{R}^d}
H(X_{\mathbf{s}}; \sigma_{s_1} u_1, \ldots,
\sigma_{s_d} u_d) \varphi_d(\mathbf{u}) \,d
\mathbf{u}.
\]
It is not essential that the weight process equals the diffusion
process $X$. Instead, we may consider any
$k$-dimensional $(\mathcal F_t)$-adapted It\^o semimartingale of type
(\ref{diffusion}). We leave
the details to the interested reader.
\end{remark}

\section{Stable central limit theorem}\label{secSCLT}
In this section we present a functional stable central limit theorem
associated with the convergence in (\ref{LLN}).
\subsection{Stable convergence}
The concept
of stable convergence of random variables was originally introduced by
Renyi \cite{R}. For properties of stable convergence, we refer to
\cite{AE,PV}. We recall the definition of stable convergence: let
$(Y_n)_{n\in\mathbb{N}}$
be a sequence of random variables defined on $(\Omega,\mathcal
{F},\mathbb{P})$ with values in a Polish space $(E, \mathcal{E})$.
We say that $Y_n$ converges stably with limit $Y$, written $Y_n
\stackrel{\mathrm{st}}{\longrightarrow}Y$, where $Y$ is defined on an extension
$(\Omega',\mathcal{F}',\mathbb{P}')$ of the original probability
space $(\Omega,\mathcal{F},\mathbb{P})$,
if and only if for any bounded, continuous function $g$ and any bounded
$\mathcal{F}$-measurable random variable $Z$ it holds that
%
%
\begin{equation}
\label{defstable} \mathbb{E}\bigl[ g(Y_n) Z\bigr] \rightarrow
\mathbb{E}'\bigl[ g(Y) Z\bigr], \qquad n \rightarrow\infty.
\end{equation}
Typically, we will deal with $E=\mathbb D([0,T], \mathbb{R})$
equipped with the Skorohod topology, or the uniform topology if the
process $Y$ is continuous. Notice that stable convergence is a stronger
mode of convergence than weak convergence. In fact, the statement
$Y_n \stackrel{\mathrm{st}}{\longrightarrow}Y$ is equivalent to the joint weak
convergence $(Y_n, Z) \stackrel{d}{\longrightarrow}(Y,Z)$ for any
$\mathcal F$-measurable random variable
$Z$; see, for example,~\cite{AE}.

\subsection{Central limit theorem}

For the stable central limit theorem we require a further structural
assumption on the volatility process $(\sigma_s)_{s\geq0}$.
We assume that $\sigma$ itself is a continuous It\^o semimartingale,
%
%
\begin{equation}
\label{assa} \sigma_t=\sigma_0 + \int
_0^t \tilde{a}_s \,ds + \int
_0^t \tilde{\sigma}_s
\,dW_s + \int_0^t
\tilde{v}_s \,dV_s,
\end{equation}
where the processes $(\tilde{a}_s)_{s\geq0}$, $(\tilde{\sigma
}_s)_{s\geq0}$,
$(\tilde{v}_s)_{s\geq0}$ are c\`adl\`ag, adapted and $V$ is a
Brownian motion independent of $W$. This type of condition
is motivated by potential applications. For instance, when $\sigma
_t=f(X_t)$ for a $C^2$-function $f$, then the It\^o formula
implies representation (\ref{assa}) with $\tilde{v}\equiv0$. In
fact, a condition of type (\ref{assa}) is nowadays a standard
assumption for proving stable central limit theorems for functionals of
high frequency data; see, for example,~\cite{BGJPS,J2}.
Moreover, we assume that the process $\sigma$ does not vanish, that is,
%
%
\begin{equation}
\label{nonvanish} \sigma_s \neq0 \qquad\mbox{for all $s\in[0,T]$}.
\end{equation}
We believe that this assumption is not essential, but dropping it would
make the following proofs considerably more involved and technical.
As in the previous subsection, the central limit theorems presented in
this paper are stable under stopping. This means,
we may assume, without loss of generality, that
%
%
\begin{equation}
\label{probound2} \mbox{The processes $a,\sigma, \sigma^{-1}, \tilde{a},
\tilde{\sigma}$ and $\tilde{v}$ are bounded in $(\omega,t)$.}
\end{equation}
We refer again to \cite{BGJPS}, Section~3, for a detailed justification
of this statement.

We need to introduce some further notation to describe the limiting process.
First, we will study the asymptotic properties of the empirical process
%
%
\begin{equation}
\label{gn} \mathbb G_n(t,x):= \frac{1}{\sqrt{n}} \sum
_{j=1}^{[nt]} \bigl( \mathbh{1}_{\{\alpha_j^n\leq x\}} -
\Phi_{\sigma_{(j-1)/{n}}}(x) \bigr),
\end{equation}
where $\alpha_j^n$ is defined at (\ref{alpha}).
This process is of crucial importance for proving the stable central
limit theorem for the $U$-statistic $U(H)_t^n$.
We start with the derivation of some useful inequalities for the
process $\mathbb G_n$.

%
\begin{lemma} \label{lem1}
For any even number $p\geq2$ and $x,y\in\mathbb{R}$, we obtain the
inequalities
%
%
\begin{eqnarray}
\label{gnineq1}\mathbb{E}\Bigl[\sup_{t\in[0,T]} \bigl|\mathbb
G_n (t, x)\bigr|^p\Bigr]&\leq& C_{T,p} \phi(x),
\\
\label{gnineq2} \mathbb{E}\Bigl[\sup_{t\in[0,T]} \bigl|\mathbb
G_n (t, x) - \mathbb G_n (t, y)\bigr|^p\Bigr]&\leq&
C_{T,p} |x-y|,
\end{eqnarray}
where $\phi\dvtx\mathbb{R}\rightarrow\mathbb{R}$ is a bounded function
(that depends on $p$ and $T$) with exponential decay at $\pm\infty$.
\end{lemma}
\begin{pf} Recall that the processes $\sigma$ and $\sigma^{-1}$ are
assumed to be bounded. We begin with inequality (\ref{gnineq1}). For
any given $x\in\mathbb{R}$, $(\mathbb G_n (t, x))_{t\in[0,T]}$
is an \mbox{$(\mathcal F_{[nt]/n})$-}martingale. Hence, the discrete
Burkh\"older inequality implies that
\[
\mathbb{E}\Bigl[\sup_{t\in[0,T]} \bigl|\mathbb G_n (t,
x)\bigr|^p\Bigr] \leq C_{T,p} \mathbb{E} \Biggl[ \Biggl| \sum
_{j=1}^{[nT]} \zeta_j^n
\Biggr|^{p/2} \Biggr]
\]
with $\zeta_j^n:= n^{-1}(\mathbh{1}_{\{\alpha_j^n\leq x\}} - \Phi
_{\sigma_{(j-1)/n}}(x) )^2$.
Recalling that $p\geq2$ is an even number und applying the H\"older
inequality, we deduce that
\begin{eqnarray*}
\Biggl| \sum_{j=1}^{[nT]} \zeta_j^n
\Biggr|^{p/2} &\leq& C_T n^{-1} \sum
_{j=1}^{[nT]} \bigl(\mathbh{1}_{\{\alpha_j^n\leq x\}} -
\Phi_{\sigma_{(j-1)/{n}}}(x)\bigr)^p
\\
&=&C_T n^{-1} \sum_{j=1}^{[nT]}
\sum_{k=0}^p {\pmatrix{p
\cr
k}}
(-1)^k \Phi_{\sigma_{(j-1)/{n}}}^k (x) \mathbh{1}_{\{\alpha
_j^n\leq x\}}.
\end{eqnarray*}
Thus we conclude that
\[
\mathbb{E}\Bigl[\sup_{t\in[0,T]} \bigl|\mathbb G_n (t,
x)\bigr|^p\Bigr] \leq C_{T,p} \sup_{s\in[0,T]}
\mathbb{E}\bigl[\Phi_{\sigma_{s}}(x) \bigl(1- \Phi_{\sigma
_{s}}(x)
\bigr)^p\bigr]=: C_{T,p} \phi(x),
\]
where the function $\phi$ obviously satisfies our requirements.
This completes the proof of (\ref{gnineq1}). By exactly the same
methods we obtain, for any $x\geq y$,
\begin{eqnarray*}
&& \mathbb{E}\Bigl[\sup_{t\in[0,T]} \bigl|\mathbb G_n (t, x) -
\mathbb G_n (t, y)\bigr|^p\Bigr]
\\
&&\qquad \leq C_{T,p} \sup_{s\in[0,T]} \mathbb{E}\bigl[\bigl(
\Phi_{\sigma_{s}}(x) - \Phi_{\sigma_{s}}(y)\bigr) \bigl(1- \bigl(
\Phi_{\sigma_{s}}(x) - \Phi_{\sigma_{s}}(y)\bigr)\bigr)^p
\bigr].
\end{eqnarray*}
Since $\sigma$ and $\sigma^{-1}$ are both bounded, there exists a
constant $M>0$ such that
\[
\sup_{s\in[0,T]}\bigl|\Phi_{\sigma_{s}}(x) -\Phi_{\sigma_{s}}(y)\bigr|
\leq|x-y| \sup_{M^{-1}\leq z\leq M, y\leq r\leq x} \varphi_z (r).
\]
This immediately gives (\ref{gnineq2}).
\end{pf}

Our next result presents a functional stable central limit theorem for
the process $\mathbb G_n$ defined at (\ref{gn}).

%
\begin{proposition} \label{prop2}
We obtain the stable convergence
\[
\mathbb G_n(t,x) \stackrel{\mathrm{st}} {\longrightarrow}\mathbb
G(t,x)
\]
on $\mathbb D([0,T])$ equipped with the uniform topology,
where the convergence is functional in $t\in[0,T]$ and in finite
distribution sense in $x\in\mathbb{R}$. The limiting process $\mathbb G$
is defined on an extension $(\Omega',\mathcal{F}',\mathbb{P}')$ of
the original probability space $(\Omega,\mathcal{F},\mathbb{P})$ and
it is Gaussian conditionally on $\mathcal F$.
Its conditional drift and covariance kernel are given by
\begin{eqnarray*}
&& \mathbb{E}'\bigl[\mathbb G (t,x)|\mathcal F\bigr] = \int
_0^t \widebar{\Phi}_{\sigma_s} (x)
\,dW_s,
\\
&& \mathbb{E}'\bigl[\mathbb G(t_1,x_1)
\mathbb G(t_2,x_2) |\mathcal F\bigr] -
\mathbb{E}'\bigl[\mathbb G(t_1,x_1) |
\mathcal F\bigr] \mathbb{E}'\bigl[\mathbb G(t_2,x_2)
|\mathcal F\bigr]
\\
&&\qquad = \int_0^{t_1 \wedge t_2} \Phi_{\sigma_s}
(x_1 \wedge x_2) - \Phi_{\sigma_s}
(x_1) \Phi_{\sigma_s} (x_2) - \widebar{\Phi
}_{\sigma_s} (x_1) \widebar{\Phi}_{\sigma_s}
(x_2)\,ds,
\end{eqnarray*}
where $\widebar{\Phi}_{z} (x) = \mathbb{E}[V\mathbh{1}_{\{zV\leq
x\}}]$ with $V\sim\mathcal N(0,1)$.
\end{proposition}
\begin{pf} Recall that due to (\ref{probound2}) the process $\sigma$
is bounded in $(\omega,t)$.
[However, note that we do not require the condition (\ref{assa}) to hold.]
For any given \mbox{$x_1, \ldots, x_k\in\mathbb{R}$}, we need to prove the
functional stable convergence
\[
\bigl(\mathbb G_n(\cdot,x_1), \ldots, \mathbb
G_n(\cdot,x_k)\bigr) \stackrel{\mathrm{st}} {\longrightarrow}
\bigl(\mathbb G(\cdot,x_1), \ldots, \mathbb G(\cdot,x_k)
\bigr).
\]
We write $\mathbb G_n(t,x_l) = \sum_{j=1}^{[nt]} \chi_{j,l}^n$ with
\[
\chi_{j,l}^n:= \frac{1}{\sqrt{n}} \bigl(
\mathbh{1}_{\{\alpha
_j^n\leq x_l\}} - \Phi_{\sigma_{(j-1)/{n}}}(x_l) \bigr), \qquad
1\leq l\leq k.
\]
According to \cite{JS}, Theorem IX.7.28, we need to show that
%
%
\begin{eqnarray}
 && \sum_{j=1}^{[nt]} \mathbb{E}
\bigl[\chi_{j,r}^n \chi_{j,l}^n |
\mathcal F_{(j-1)/{n}}\bigr]
\nonumber\\[-8pt]\label{p1} \\[-8pt]
&&\qquad \stackrel{\mathbb{P}} {\longrightarrow} \int
_0^t \bigl( \Phi_{\sigma_s}
(x_r \wedge x_l) - \Phi_{\sigma_s}
(x_r) \Phi_{\sigma_s} (x_l) \bigr) \,ds,\nonumber
\\
 \label{p2} && \sum_{j=1}^{[nt]} \mathbb{E}
\bigl[\chi_{j,l}^n \Delta_j^n W |
\mathcal F_{(j-1)/{n}}\bigr] \stackrel{\mathbb{P}} {\longrightarrow
} \int
_0^t \widebar{\Phi}_{\sigma_s}
(x_l) \,ds,
\\
\label{p3} && \sum_{j=1}^{[nt]} \mathbb{E}
\bigl[\bigl|\chi_{j,l}^n\bigr|^2 1_{\{
|\chi_{j,l}^n|>\varepsilon\}} |
\mathcal F_{(j-1)/{n}}\bigr] \stackrel{\mathbb{P}} {\longrightarrow
}0\qquad
\mbox{for all $\varepsilon>0$,}
\\
\label{p4} && \sum_{j=1}^{[nt]} \mathbb{E}
\bigl[\chi_{j,l}^n \Delta_j^n N |
\mathcal F_{(j-1)/{n}}\bigr] \stackrel{\mathbb{P}} {\longrightarrow}0,
\end{eqnarray}
where $1\leq r,l\leq d$ and the last condition must hold for all
bounded continuous martingales $N$ with $[W,N]=0$.
The convergence in (\ref{p1}) and (\ref{p2}) is obvious, since
$\Delta_j^n W$ is independent of $\sigma_{(j-1)/n}$.
We also have that
\[
\sum_{j=1}^{[nt]} \mathbb{E}\bigl[\bigl|
\chi_{j,l}^n\bigr|^2 1_{\{|\chi
_{j,l}^n|>\varepsilon\}} | \mathcal
F_{(j-1)/{n}}\bigr] \leq\varepsilon^{-2} \sum
_{j=1}^{[nt]} \mathbb{E}\bigl[\bigl|\chi_{j,l}^n\bigr|^4
| \mathcal F_{(j-1)/{n}}\bigr] \leq Cn^{-1},
\]
which implies (\ref{p3}). Finally, let us prove (\ref{p4}). We fix
$l$ and define
$M_u:= \mathbb{E}[\chi_{j,l}^n | \mathcal F_{u}]$ for $u\geq
(j-1)/n$. By the martingale representation theorem we deduce the
identity
\[
M_u = M_{(j-1)/{n}} + \int_{(j-1)/{n}}^u
\eta_s \,dW_s
\]
for a suitable predictable process $\eta$. By the It\^o isometry we
conclude that
\[
\mathbb{E}\bigl[\chi_{j,l}^n \Delta_j^n
N | \mathcal F_{(j-1)/{n}}\bigr] = \mathbb{E}\bigl[M_{{j}/{n}}
\Delta_j^n N | \mathcal F_{(j-1)/{n}}\bigr] =
\mathbb{E}\bigl[ \Delta_j^n M \Delta_j^n
N | \mathcal F_{(j-1)/{n}}\bigr] =0.
\]
This completes the proof of Proposition~\ref{prop2}.
\end{pf}

We suspect that the stable convergence in Proposition~\ref{prop2} also
holds in the functional sense in the $x$ variable.
However, proving tightness (even on compact sets) turns out to be a
difficult task. In particular,
inequality (\ref{gnineq2}) is not sufficient for showing tightness.
%
%
\begin{remark} \label{rem2}
We highlight some probabilistic properties of the limiting process
$\mathbb G$ defined in Proposition~\ref{prop2}.
\begin{longlist}[(iii)]
\item[(i)] Proposition~\ref{prop2} can be reformulated as follows.
Let $x_1, \ldots, x_k\in\mathbb{R}$ be arbitrary real numbers.
Then it holds that
\[
\bigl(\mathbb G_n(\cdot,x_1), \ldots, \mathbb
G_n(\cdot,x_k)\bigr) \stackrel{\mathrm{st}} {\longrightarrow}\int
_0^{\cdot} v_s \,dW_s+ \int
_0^{\cdot} w_s^{1/2}
\,dW'_s,
\]
where $W'$ is a $k$-dimensional Brownian motion independent of
$\mathcal F$, and $v$ and $w$ are $\mathbb{R}^k$-valued
and $\mathbb{R}^{k \times k}$-valued processes, respectively, with coordinates
\begin{eqnarray*}
v_s^r &=& \widebar{\Phi}_{\sigma_s}
(x_r),
\\
w_s^{rl} &=& \Phi_{\sigma_s} (x_r \wedge
x_l) - \Phi_{\sigma_s} (x_r) \Phi_{\sigma_s}
(x_l) - \widebar{\Phi}_{\sigma_s} (x_r)
\widebar{\Phi}_{\sigma_s} (x_l),
\end{eqnarray*}
for $1\leq r,l\leq k$. This type of formulation appears in \cite{JS},
Theorem IX.7.28. In particular,
$(\mathbb G(\cdot,x_l))_{1\leq l\leq k}$ is a $k$-dimensional
martingale.

\item[(ii)] It is obvious from (i) that $\mathbb G$ is continuous in
$t$. Moreover, $\mathbb G$
is also continuous in $x$. This follows from Kolmogorov's criterion and
the inequality ($y\leq x$)
\begin{eqnarray*}
&& \mathbb{E}'\bigl[\bigl|\mathbb G(t,x) -\mathbb{G}(t,y)\bigr|^p
\bigr]
\\
&&\qquad \leq C_p \mathbb{E} \biggl[ \biggl( \int_0^t
\bigl\{ \Phi_{\sigma_s}(x) - \Phi_{\sigma_s}(y) - \bigl(
\Phi_{\sigma_s}(x) - \Phi_{\sigma_s}(y)\bigr)^2 \bigr\} \,ds
\biggr)^{p/2} \biggr]
\\
&&\qquad \leq C_p (x-y)^{p/2},
\end{eqnarray*}
for any $p>0$, which follows by the Burkh\"older inequality. In
particular, $\mathbb G(t,\cdot)$
has H\"older continuous paths of order $1/2-\varepsilon$, for any
$\varepsilon\in(0, 1/2)$.

\item[(iii)] A straightforward computation [cf.~(\ref{gnineq1})]
shows that the function
$\mathbb{E}[ \sup_{t\in[0,T]}\mathbb G(t,x)^2]$ has
exponential decay as $x\rightarrow\pm\infty$.
Hence, for any function $f \in C^1_p(\mathbb{R})$, we have
\[
\int_{\mathbb{R}} f(x) \mathbb G(t,dx) < \infty, \qquad\mbox{a.s.}
\]
If $f$ is an even function, we also have that
\[
\int_{\mathbb{R}} f(x) \mathbb G(t,dx) = \int_{\mathbb{R}}
f(x) \bigl( \mathbb G(t,dx) - \mathbb{E}'\bigl[\mathbb G(t,dx) |
\mathcal{F}\bigr]\bigr),
\]
since
\[
\int_{\mathbb{R}} f(x) \mathbb{E}'\bigl[\mathbb
G(t,dx)|\mathcal{F}\bigr] = \int_0^t \biggl(
\int_{\mathbb{R}} f(x) \widebar{\Phi}_{\sigma_s} (dx) \biggr)
\,dW_s,
\]
and, for any $z>0$,
\[
\int_{\mathbb{R}} f(x) \widebar{\Phi}_{z} (dx) = \int
_{\mathbb
{R}} xf(x) \varphi_z(x) \,dx =0,
\]
because $f \varphi_z $ is an even function. The same argument applies
for $z<0$. Furthermore,
the integration by parts formula and the aforementioned argument imply
the identity
\begin{eqnarray*}
&& \mathbb{E}' \biggl[ \biggl|\int_{\mathbb{R}} f(x) \mathbb
G(t,dx) \biggr|^2 \Big| \mathcal F \biggr]
\\
&&\qquad = \int_0^{t} \biggl( \int_{\mathbb{R}^2}
f'(x) f'(y) \bigl( \Phi_{\sigma_s} (x \wedge y) -
\Phi_{\sigma_s} (x) \Phi_{\sigma_s} (y) \bigr) \,dx\,dy \biggr) \,ds.
\end{eqnarray*}
We remark that, for any $z \neq0$, we have
\[
\operatorname{var}\bigl[f (V)\bigr]=\int_{\mathbb{R}^2} f'(x)
f'(y) \bigl( \Phi_{z} (x \wedge y) -
\Phi_{z} (x) \Phi_{z} (y) \bigr) \,dx\,dy
\]
with $V\sim\mathcal N(0,z^2)$.
\end{longlist}
\end{remark}
Now, we present a functional stable central limit theorem of the
$U$-statistic $U^{\prime n}_t (H)$ given at (\ref{uprime}), which is
based on the approximative
quantities $(\alpha_j^n)_{1\leq j\leq[nt]}$ defined at
(\ref{alpha}).

%
\begin{proposition} \label{prop3}
Assume that conditions (\ref{assa}), (\ref{nonvanish}) and (\ref
{probound2}) hold.
Let $H\in C_p^1(\mathbb{R}^d)$ be a symmetric function that is even in
each (or, equivalently, in one) argument.
Then we obtain the functional stable convergence
%
%
\begin{equation}
\label{uprimeclt} \sqrt{n} \bigl(U^{\prime n} (H) - U(H)\bigr) \stackrel
{\mathrm{st}} {
\longrightarrow}L,
\end{equation}
where
%
%
\begin{equation}
\label{Lt} L_t=d\int_{\mathbb{R}^d} H(x_1,
\ldots, x_d) \mathbb G(t,dx_1) F(t,dx_2)
\cdots F(t,dx_d).
\end{equation}
The convergence takes place in $\mathbb D([0,T])$ equipped with the
uniform topology. Furthermore, $\mathbb G$ can be replaced by
$\mathbb G - \mathbb{E}'[\mathbb G|\mathcal F]$ without changing the
limit and, consequently, $L$ is a centered Gaussian process,
conditionally on $\mathcal F$.
\end{proposition}
\begin{pf} First of all, we remark that
\[
\int_{\mathbb{R}} H(x_1, \ldots, x_d)
\mathbb{E}'\bigl[\mathbb G(t,dx_1) | \mathcal F\bigr]
=0
\]
follows from Remark~\ref{rem2}(iii). The main part of the proof is
divided into five steps:
\begin{longlist}[(iii)]
\item[(i)] In Section~\ref{subsecsub33} we will show that
under condition (\ref{assa}) we have
%
%
\begin{equation}
\label{riemannphi} \sqrt{n} \biggl( U(H)_{t} - \int_{\mathbb{R}^d}
H(\mathbf x) \widebar{F}_n^{\otimes d} (t, d\mathbf x) \biggr)
\stackrel{\mathrm{u.c.p.}} {\longrightarrow}0
\end{equation}
with
\[
\widebar{F}_n (t, x):= \frac{1}n \sum
_{j=1}^{[nt]} \Phi_{\sigma
_{(j-1)/{n}}} (x).
\]
Thus, we need to prove the stable convergence $L^n \stackrel
{\mathrm{st}}{\longrightarrow}L$ for
%
%
\begin{equation}
\label{ltn} L_t^n:= \sqrt{n} \biggl(U^{\prime n}_{t}
(H) - \int_{\mathbb{R}^d} H(\mathbf x) \widebar{F}_n^{\otimes d}
(t, d\mathbf x) \biggr).
\end{equation}
Assume that the function $H\in C^1(\mathbb{R}^d)$ has compact support.
Recalling definition~(\ref{gn}) of the empirical process $\mathbb
G_n$, we obtain the identity
\[
L_t^n = \sum_{l=1}^d
\int_{\mathbb{R}^d} H(\mathbf x) \mathbb G_n (t,
dx_l) \prod_{m=1}^{l-1}
F_n(t, dx_m) \prod_{m=l+1}^{d}
\widebar{F}_n(t, dx_m).
\]
In step (iv) we will show that both $F_n(t, dx_m)$ and $\widebar
{F}_n(t, dx_m)$ can be replaced by $F(t, dx_m)$
without affecting the limit. In other words, $L^n - L^{\prime n}
\stackrel{\mathrm{u.c.p.}}{\longrightarrow}
0$ with
\[
L_t^{\prime n}:= \sum_{l=1}^d
\int_{\mathbb{R}^d} H(\mathbf x) \mathbb G_n (t,
dx_l) \prod_{m\neq l} F(t,
dx_m).
\]
But, since $H$ is symmetric, we readily deduce that
\[
L_t^{\prime n} = d \int_{\mathbb{R}^d} H(\mathbf x)
\mathbb G_n (t, dx_1) \prod_{m = 2}^d
F(t, dx_m).
\]
The random measure $F(t, x)$ has a Lebesgue density in $x$
due to assumption (\ref{nonvanish}), which we denote by $F'(t, x)$.
The integration by parts formula implies that
\[
L_t^{\prime n} = -d \int_{\mathbb{R}^{d}}
\partial_1 H(\mathbf x) \mathbb G_n (t, x_1)
\prod_{m = 2}^d F'(t,
x_m) \,d\mathbf x,
\]
where $\partial_l H$ denotes the partial derivative of $H$ with
respect to $x_l$. This identity completes step (i).

\item[(ii)] In this step we will start proving the stable convergence
$L^{\prime n} \stackrel{\mathrm{st}}{\longrightarrow}L$ [the function $H\in
C^1(\mathbb{R}^d)$
is still assumed to have compact support]. Since the stable convergence
$\mathbb G_n \stackrel{\mathrm{st}}{\longrightarrow}
\mathbb G$ does not hold in the functional sense in the $x$ variable,
we need to overcome
this problem by a Riemann sum approximation. Let the support of $H$ be
contained in $[-k,k]^d$. Let $-k=z_0<\cdots<z_l=k$ be the equidistant
partition of the interval $[-k,k]$.
We set
\[
Q (t,x_1):= \int_{\mathbb{R}^{d-1}} \partial_1
H(x_1, \ldots,x_d) \prod_{m = 2}^d
F'(t, x_m) \,dx_2 \cdots dx_d,
\]
and define the approximation
of $L_t^{\prime n}$ via
\[
L_t^{\prime n} (l) = - \frac{2d k}{l} \sum
_{j=0}^l Q (t,z_j) \mathbb
G_n (t, z_{j}).
\]
Proposition~\ref{prop2} and the properties of stable convergence imply that
\[
\bigl(Q (\cdot,z_j), \mathbb G_n (\cdot,
z_{j}) \bigr)_{0\leq j\leq
l} \stackrel{\mathrm{st}} {\longrightarrow} \bigl(Q
(\cdot,z_j), \mathbb G (\cdot, z_{j})
\bigr)_{0\leq j\leq l}.
\]
Hence, we deduce the stable convergence
\[
L_{\cdot}^{\prime n} (l) \stackrel{\mathrm{st}} {\longrightarrow}L_{\cdot
}(l):=
- \frac{2d k}{l} \sum_{j=0}^l Q (
\cdot,z_j) \mathbb G (\cdot, z_{j})
\]
as $n\rightarrow\infty$, for any fixed $l$. Furthermore, we obtain
the convergence
\[
L(l) \stackrel{\mathrm{u.c.p.}} {\longrightarrow}L
\]
as $l\rightarrow\infty$, where we reversed all above transformations.
This convergence completes step (ii).

\item[(iii)] To complete the proof of the stable convergence $L^{\prime n}
\stackrel{\mathrm{st}}{\longrightarrow}L$,
we need to show that
\[
\lim_{l\rightarrow\infty} \limsup_{n\rightarrow\infty} \sup
_{t\in[0,T]}\bigl|L_{t}^{\prime n} (l)- L_{t}^{\prime n}
\bigr| =0,
\]
where the limits are taken in probability. With $h=l/2k$ we obtain that
\[
\bigl|L_{t}^{\prime n} (l)- L_{t}^{\prime n} \bigr| = d
\biggl\llvert\int_{\mathbb
{R}} \bigl\{ Q \bigl(t,[xh]/h\bigr)
\mathbb G_n \bigl(t, [xh]/h\bigr) - Q (t,x) \mathbb G_n
(t, x) \bigr\} \,dx\biggr\rrvert.
\]
Observe that
%
%
\begin{equation}
\label{fprimeapp} \sup_{t\in[0,T]} \bigl|F'(t,
x_m)\bigr| = \int_{0}^T
\varphi_{\sigma_s}(x_m) \,ds \leq T \sup_{M^{-1}\leq z\leq M}
\varphi_{z}(x_m),
\end{equation}
where $M$ is a positive constant with $M^{-1}\leq|\sigma|\leq M$.
Recalling the definition of $Q (t,x)$ we obtain that
%
%
\begin{eqnarray}
\label{qest}
\sup_{t\in[0,T]} \bigl|Q (t,x)\bigr| &\leq& C_T,
\nonumber\\[-8pt]\\[-8pt]
\sup_{t\in[0,T]} \bigl|Q (t,x) - Q \bigl(t,[xh]/h\bigr)\bigr| &\leq&
C_T \eta\bigl(h^{-1}\bigr),\nonumber
\end{eqnarray}
where $\eta(\varepsilon):= \sup\{|\partial_1 H(\mathbf y_1) -
\partial_1 H(\mathbf y_2)|\dvtx
\|\mathbf y_1-\mathbf y_2\|\leq\varepsilon, \mathbf y_1,\mathbf
y_2\in[-k,k]^d\}$ denotes the modulus of continuity
of the function $\partial_1 H$. We also deduce by Lemma~\ref{lem1} that
%
%
\begin{eqnarray}
\label{gnest1} \mathbb{E}\Bigl[\sup_{t\in[0,T]} \bigl|\mathbb
G_n (t, x)\bigr|^p\Bigr]&\leq& C_T,
\\
\label{gnest2} \mathbb{E}\Bigl[\sup_{t\in[0,T]} \bigl|\mathbb
G_n (t, x) - \mathbb G_n\bigl(t,[xh]/h
\bigr)\bigr|^p\Bigr] &\leq& C_T h^{-1},
\end{eqnarray}
for any even number $p\geq2$.
Combining inequalities (\ref{qest}), (\ref{gnest1}) and (\ref
{gnest2}), we deduce the convergence
\[
\lim_{l\rightarrow\infty} \limsup_{n\rightarrow\infty} \mathbb{E}\Bigl[
\sup_{t\in[0,T]}\bigl|L_{t}^{\prime n} (l)-
L_{t}^{\prime n} \bigr|\Bigr] =0
\]
using that $Q(t,\cdot)$ has compact support contained in $[-k,k]$.
Hence, $L^{\prime n} \stackrel{\mathrm{st}}{\longrightarrow}L$, and we are
done.

\item[(iv)] In this step we will prove the convergence
\[
L^n - L^{\prime n} \stackrel{\mathrm{u.c.p.}} {
\longrightarrow}0.
\]
This difference can be decomposed into several terms; in the following
we will treat a typical representative
(all other terms are treated in exactly the same manner).
For $l<d$ define
\begin{eqnarray*}
R_t^n (l)&:=&\int_{\mathbb{R}^d} H(\mathbf x)
\mathbb G_n (t, dx_l) \prod_{m=1}^{l-1}
F_n(t, dx_m)
\\
&&{}\times \prod_{m=l+1}^{d-1}
\widebar{F}_n(t, dx_m) \bigl[\widebar{F}_n(t,
dx_d) - F(t, dx_d)\bigr].
\end{eqnarray*}
Now, we use the integration by parts formula to obtain that
\[
R_t^n (l)= \int_{\mathbb{R}}
N_n(t, x_l) \mathbb G_n (t,
x_l) \,dx_l,
\]
where
\begin{eqnarray*}
N_n(t, x_l) &=& \int_{\mathbb{R}^{d-1}}
\partial_l H(\mathbf x) \prod_{m=1}^{l-1}
F_n(t, dx_m)
\\
&&{}\times \prod_{m=l+1}^{d-1}
\widebar{F}_n(t, dx_m) \bigl[\widebar{F}_n(t,
dx_d) - F(t, dx_d)\bigr].
\end{eqnarray*}
As in step (iii) we deduce for any even $p\geq2$,
\begin{eqnarray*}
\mathbb{E}\Bigl[\sup_{t\in[0,T]} \bigl|\mathbb G_n (t,
x_l)\bigr|^p\Bigr]&\leq& C_p,
\\
\mathbb{E}
\Bigl[\sup_{t\in[0,T]} \bigl|N_n(t, x_l)\bigr|^p
\Bigr]&\leq& C_p.
\end{eqnarray*}
Recalling that the function $H$ has compact support and applying the
dominated convergence theorem,
it is sufficient to show that
\[
N_n(\cdot, x_l) \stackrel{\mathrm{u.c.p.}} {
\longrightarrow}0,
\]
for any fixed $x_l$. But this follows immediately from Lemma~\ref
{vague}, since
\[
F_n(\cdot, x) \stackrel{\mathrm{u.c.p.}} {\longrightarrow}F(\cdot,
x), \qquad\widebar{F}_n(\cdot, x) \stackrel{\mathrm{u.c.p.}} {
\longrightarrow}F(\cdot, x),
\]
for any fixed $x\in\mathbb{R}$, and $\partial_l H$ is a continuous
function with compact support. This finishes the proof
of step (iv).

\item[(v)] Finally, let $H \in C_p^1(\mathbb{R}^d)$ be arbitrary. For any
$k\in\mathbb{N}$, let $H_k \in C_p^1(\mathbb{R}^d)$ be a function
with $H_k=H$ on
$[-k,k]^d$ and $H_k=0$ on $([-k-1,k+1]^d)^c$. Let us denote by $L_t^n
(H)$ and $L_t (H)$ the processes
defined by (\ref{ltn}) and (\ref{Lt}), respectively, that are
associated with a given function $H$. We know from the previous steps
that
\[
L^n (H_k) \stackrel{\mathrm{st}} {\longrightarrow}L(H_k)
\]
as $n\rightarrow\infty$, and $L(H_k) \stackrel{\mathrm
{u.c.p.}}{\longrightarrow}L(H)$ as $k\rightarrow
\infty$.
So, we are left to proving that
\[
\lim_{k\rightarrow\infty} \limsup_{n\rightarrow\infty} \sup
_{t\in[0,T]}\bigl|L^n_t (H_k)-
L^n_t(H)\bigr| =0,
\]
where the limits are taken in probability. As in steps (ii) and (iii)
we obtain the identity
\begin{eqnarray*}
&& L^n_t (H_k)- L^n_t(H)
\\
&&\qquad = \sum_{l=1}^d \int_{\mathbb{R}^d}
\partial_l (H-H_k) (\mathbf x) \mathbb G_n
(t, x_l) \,dx_l \prod_{m=1}^{l-1}
F_n(t, dx_m) \prod_{m=l+1}^{d}
\widebar{F}_n(t, dx_m)
\\
&&\qquad =: \sum_{l=1}^d Q^{l}
(k)_t^n.
\end{eqnarray*}
We deduce the inequality
\begin{eqnarray*}
\bigl|Q^{l} (k)_t^n\bigr| &\leq& n^{-(l-1)} \sum
_{i_1,\ldots,i_{l-1}=1}^{[nt]} \int_{\mathbb{R}^{d-l+1}}
\bigl|\partial_l (H-H_k) \bigl(\alpha_{i_1}^n,
\ldots, \alpha_{i_{l-1}}^n, x_l,\ldots,
x_d\bigr)\bigr|
\\
&&\hspace*{108pt}{}\times\bigl|\mathbb G_n (t, x_l)\bigr| \prod
_{m=l+1}^{d} \widebar{F}'_{n}(t,
x_m) \,dx_l\cdots dx_d.
\end{eqnarray*}
We remark that $\partial_l(H_k-H)$ vanishes if all arguments lie in
the interval $[-k,k]$. Hence
\begin{eqnarray*}
\bigl|Q^{l} (k)_t^n\bigr| &\leq& n^{-(l-1)} \sum
_{i_1,\ldots,i_{l-1}=1}^{[nt]} \int_{\mathbb{R}^{d-l+1}} \bigl|
\partial_l (H-H_k) \bigl(\alpha_{i_1}^n,
\ldots, \alpha_{i_{l-1}}^n, x_l,\ldots,
x_d\bigr)\bigr|
\\
&&\hspace*{109pt}{}\times\Biggl(\sum_{m=1}^{l-1}
\mathbh{1}_{\{|\alpha_{i_m}^n|>k\}} + \sum_{m=l}^{d}
\mathbh{1}_{\{|x_m|>k\}} \Biggr)
\\
&&\hspace*{109pt}{}\times\bigl|\mathbb G_n (t, x_l)\bigr|
\prod_{m=l+1}^{d} \widebar{F}'_{n}(t,
x_m) \,dx_l\cdots dx_d.
\end{eqnarray*}
Now, applying Lemma~\ref{lem1}, (\ref{pmoments}), (\ref{fprimeapp})
and the Cauchy--Schwarz inequality, we deduce that
\begin{eqnarray*}
&& \mathbb{E}\Bigl[\sup_{t\in[0,T]}\bigl|Q^{l}
(k)_t^n\bigr|\Bigr]
\\
&&\qquad \leq  C_T \int_{\mathbb{R}^{d-l+1}} \Biggl((l-1)\sup
_{M^{-1}\leq z\leq M} \bigl(1-\Phi_z(k)\bigr) + \sum
_{m=l}^{d} \mathbh{1}_{\{|x_m|>k\}}
\Biggr)^{1/2}
\\
&&\hspace*{78pt} {}\times\psi(x_l,\ldots, x_d)
\phi(x_l) \prod_{m=l+1}^{d} \sup
_{M^{-1} \leq z\leq M} \varphi_z(x_m)
\,dx_l\cdots dx_d,
\end{eqnarray*}
for some bounded function $\phi$ with exponential decay at $\pm\infty
$ and a function $\psi\in C_p^0 (\mathbb{R}^{d-l+1})$.
Hence
\[
\int_{\mathbb{R}^{d-l+1}} \psi(x_l,\ldots, x_d)
\phi(x_l) \prod_{m=l+1}^{d} \sup
_{M^{-1} \leq z\leq M} \varphi_z(x_m)
\,dx_l\cdots dx_d< \infty,
\]
and we conclude that
\[
\lim_{k\rightarrow\infty} \limsup_{n\rightarrow\infty} \mathbb{E}\Bigl[
\sup_{t\in[0,T]}\bigl|Q^{l} (k)_t^n\bigr|
\Bigr] =0.
\]
This finishes step (v), and we are done with the proof of Proposition~\ref{prop3}.\quad\qed
\end{longlist}\noqed
\end{pf}

Notice that an additional $\mathcal{F}$-conditional bias would appear
in the limiting process $L$ if we would drop
the assumption that $H$ is even in each coordinate. The corresponding
asymptotic theory for the case $d=1$
has been studied in \cite{KP}; see also \cite{J3}.

%
\begin{remark}
Combining limit theorems for semimartingales with the empirical
distribution function approach is probably the most efficient way
of proving Proposition~\ref{prop3}. Nevertheless, we shortly comment
on alternative methods of proof.

Treating the multiple sum in the definition of $U^{\prime n} (H)$
directly is relatively complicated, since at a certain
stage of the proof one will have to deal with partial sums of functions
of $\alpha_j^n$ weighted by an anticipative process.
This anticipation of the weight process makes it impossible to apply
martingale methods directly.

Another approach to proving Proposition~\ref{prop3} is a \textit{pseudo}
Hoeffding decomposition. This method relies
on the application of the classical Hoeffding decomposition to
$U^{\prime n} (H)$ by pretending that
the scaling components $\sigma_{(\mathbf{i}-1)/n}$ are nonrandom.
However, since the random variables $\alpha_j^n$
are not independent when the process $\sigma$ is stochastic, the
treatment of the error term connected with the pseudo
Hoeffding decomposition will not be easy, because the usual
orthogonality arguments of the Hoeffding method
do not apply in our setting.
\end{remark}

%
\begin{remark}
In the context of Proposition~\ref{prop3} we would like to mention a
very recent work by Beutner and Z\"ahle \cite{BZ}.
They study the empirical distribution function approach to $U$- and
$V$-statistics for unbounded kernels $H$ in the classical
i.i.d. or weakly dependent setting. Their method relies on the
application of the functional delta method
for quasi-Hadamard differentiable functionals. In our setting it would
require the functional convergence
\[
\mathbb G_n(t, \cdot) \stackrel{\mathrm{st}} {\longrightarrow}\mathbb G(t,
\cdot),
\]
where the convergence takes place in the space of c\`adl\`ag functions
equipped with the weighted sup-norm
$\|f\|_\lambda:= \sup_{x\in\mathbb{R}} |(1+|x|^\lambda) f(x)|$ for
some $\lambda>0$.
Although we do not really require such a strong result in our framework
(as can be seen from the proof of Proposition
\ref{prop3}), it would be interesting to prove this type of
convergence for functionals of high frequency data; cf. the comment
before Remark~\ref{rem2}.
\end{remark}
To conclude this section, we finally present the main result: A
functional stable central limit theorem for the original
$U$-statistic $U(H)^n$.

%
\begin{theorem} \label{th2}
Assume that the symmetric function $H\in C_p^1(\mathbb{R}^d)$
is even in each (or, equivalently, in one) argument. If $\sigma$
satisfies conditions (\ref{assa}) and (\ref{nonvanish}), we obtain
the functional stable central limit theorem
%
%
\begin{equation}
\label{CLT} \sqrt{n} \bigl( U(H)^n - U(H) \bigr)\stackrel{\mathrm{st}} {
\longrightarrow}L,
\end{equation}
where the convergence takes place in $\mathbb D([0,T])$ equipped with
the uniform topology and the limiting
process $L$ is defined at (\ref{Lt}).
\end{theorem}
\begin{pf} In Section~\ref{subsecsub32} we will show the following
statement: under condition (\ref{assa})
it holds that
%
%
\begin{equation}
\label{sec3b} \sqrt{n}\bigl|U(H)^n - \widetilde{U}(H)^n\bigr|
\stackrel{\mathrm{u.c.p.}} {\longrightarrow}0.
\end{equation}
In view of Proposition~\ref{prop3}, it remains to prove that $\sqrt
{n}|\widetilde{U}(H)_t^n - U_t^{\prime n}(H)| \stackrel{\mathrm
{u.c.p.}}{\longrightarrow}0$.
But due to the symmetry of $H$, we obtain as in the proof of Theorem
\ref{th1}
\[
\mathbb{E}\Bigl[\sup_{t\in[0,T]}\bigl|\widetilde{U}(H)_t^n
- U_t^{\prime
n}(H)\bigr|\Bigr] \leq\frac{C_T}{n}.
\]
This completes the proof of Theorem~\ref{th2}.
\end{pf}

We remark that the stable convergence at (\ref{CLT}) is not \textit
{feasible} in its present form, since the distribution
of the limiting process $L$ is unknown. In the next section
we will explain how to obtain a feasible central limit theorem that
opens the door to statistical applications.

\section{Estimation of the conditional variance}\label{seccondvar}
In this section we pre\-sent a standard central limit theorem for the
$U$-statistic $U(H)_t^n$.
We will confine ourselves to the presentation of a result in finite
distributional sense.
According to Remark~\ref{rem2}(iii) applied to
\[
f_t(x):=d\int_{\mathbb{R}^{d-1}} H(x,x_2, \ldots,
x_d) F(t,dx_2)\cdots F(t,dx_d),
\]
the conditional variance of the limit $L_t$ is given by
\[
V_t:= \mathbb{E}'\bigl[|L_t|^2|
\mathcal F\bigr] = \int_0^t \biggl( \int
_{\mathbb{R}} f_t^2(x) \varphi_{\sigma_s}
(x) \,dx - \biggl(\int_{\mathbb{R}} f_t(x)
\varphi_{\sigma_s} (x) \,dx \biggr)^2 \biggr) \,ds.
\]
Hence, the random variable $L_t$ is nondegenerate when
\[
\operatorname{var} \bigl( \mathbb{E}\bigl[H(x_1
U_1, \ldots, x_d U_d)| U_1
\bigr] \bigr)>0, \qquad(U_1,\ldots, U_d)\sim\mathcal
N_d(0, \mathbf{I}_d),
\]
for all $x_1,\ldots, x_d\in\{\sigma_s| s\in A\subseteq[0,t]\}$ and
some set $A$ with positive Lebesgue measure.
This essentially coincides with the classical nonde\-gen\-er\-a\-cy
condition for $U$-statistics of independent random variables.

We define the functions $G_1\dvtx\mathbb{R}^{2d-1}\to\mathbb{R}$ and
$G_2\dvtx\mathbb{R}^2\times\mathbb{R}^{2d-2}\to\mathbb{R}$ by
%
%
\begin{eqnarray}
G_1(\mathbf{x}) &=& H(x_1,x_2,
\ldots,x_{d})H(x_1,x_{d+1},\ldots,x_{2d-1}),\label{eqG1}
\\
G_2(\mathbf{x};\mathbf{y}) &=& H(x_1,y_1,
\ldots,y_{d-1})H(x_2,y_d,\ldots,y_{2d-2}),\label{eqG2}
\end{eqnarray}
respectively. Then $V_t$ can be written as
\begin{eqnarray*}
\label{asyvar} V_t &=& d^2\int_{[0,t]^{2d-1}}
\rho_{\sigma_{\mathbf
{s}}}(G_1) \,d\mathbf{s}
\\
&&{} - d^2\int_{[0,t]^{2d-2}}\int_0^t
\!\int_{\mathbb{R}}\int_{\mathbb{R}}\rho_{\sigma_{\mathbf{s}}}
\bigl(G_2(x_1,x_2;\cdot)\bigr)
\varphi_{\sigma_q}(x_1)\varphi_{\sigma_q}(x_2)\,dx_1
\,dx_2 \,dq \,d\mathbf{s}.
\nonumber
\end{eqnarray*}
We denote the first and second summand on the right-hand side of the
preceding equation by $V_{1,t}$ and $V_{2,t}$, respectively. Let
$\widetilde{G}_1$ denote the symmetrization of the function $G_1$. By
Theorem~\ref{th1} it holds that
\[
V_{1,t}^n=d^2U(\widetilde{G}_1)_t^n
\stackrel{\mathrm{u.c.p.}} {\longrightarrow}d^2 U(
\widetilde{G}_1)_t = V_{1,t}.
\]
The multiple integral $V_{2,t}$ is almost in the form of the limit in
Theorem~\ref{th1}, and it is indeed possible to estimate it by a
slightly modified $U$-statistic as the following proposition shows. The
statistic presented in the following proposition is
a generalization of the bipower concept discussed, for example, in
\cite{BGJPS} in the case $d=1$.

%
\begin{proposition}
Assume that $H\in C_p^0 (\mathbb{R}^d)$. Let
\begin{eqnarray*}
V_{2,t}^n&:=&\frac{d^2}{n}\pmatrix{n
\cr
2d-2}^{-1}
\\
&&{}\times\sum_{\mathbf{i}
\in\mathcal{A}_t^n(2d-2)}\sum_{j=1}^{[nt]-1}\widetilde
{G}_2\bigl(\sqrt{n}\Delta_j^n X,\sqrt{n}
\Delta_{j+1}^n X; \sqrt{n}\Delta_{i_1}^n
X,\ldots,\sqrt{n}\Delta_{i_{2d-2}}^n X\bigr),
\end{eqnarray*}
where $\widetilde{G}_2$ denotes the symmetrization of $G_2$ with
respect to the $\mathbf{y}$-values, that is,
\[
\widetilde{G}_2(\mathbf{x};\mathbf{y}) = \frac{1}{(2d-2)! }\sum
_{\pi} G_2(\mathbf{x};\pi\mathbf{y}),
\]
for $\mathbf{x} \in\mathbb{R}^2$, $\mathbf{y} \in\mathbb
{R}^{2d-2}$, and where the sum runs over all permutations of $\{1,\ldots,2d-2\}$.
Then
\[
V_2^n \stackrel{\mathrm{u.c.p.}} {\longrightarrow}V_2.
\]
\end{proposition}
\begin{pf}
The result can be shown using essentially the same arguments as in the
proofs of Proposition~\ref{prop1} and Theorem~\ref{th1}. We provide a
sketch of the proof. Similarly to (\ref{tildeustatistic}) we define
\[
\widetilde{V}_{2,t}^n:= \frac{d^2}{n}\pmatrix{n
\cr
2d-2}^{-1}\sum_{\mathbf{i} \in\mathcal{A}_t^n(2d-2)} \sum
_{j=1}^{[nt]-1}\widetilde{G}_2\bigl(
\alpha_j^n, \alpha_{j+1}^{\prime n};
\alpha_{i_1}^n,\ldots,\alpha_{i_{2d-2}}^n
\bigr),
\]
where $\alpha_{j+1}^{\prime n}:= \sqrt{n} \sigma_{(j-1)/{n}}
\Delta_{i+1}^n W$. Analogously to (\ref{uprime}) we introduce
the random process
\[
V^{\prime n}_{2,t}:= d^2\int_{\mathbb{R}^{2d-2}}
\int_{\mathbb{R}^2} \widetilde{G}_2(\mathbf{x};\mathbf{y})
\widetilde{F}_n(t,d\mathbf{x})F_n^{\otimes(2d-2)}(t,d
\mathbf{y}),
\]
where
\[
\widetilde{F}_n(t,x_1,x_2) =
\frac{1}{n}\sum_{j=1}^{[nt]-1} \mathbh
{1}_{\{\alpha_j^n \le x_1\}}\mathbh{1}_{\{\alpha_{j+1}^{\prime n}
\le x_2\}}.
\]
Writing out $V^{\prime n}_{2,t}$ as a multiple sum over nondecreasing
multi-indices in the $\mathbf{y}$ arguments, one observes as before
that $V^{\prime n}_{2,t}$ and $\widetilde{V}{}^n_{2,t}$ differ in at
most $O(n^{2d-3})$ summands. Therefore, using the same argument as in
the proof of Theorem~\ref{th1}
\[
\widetilde{V}{}^n_{2,t} - V^{\prime n}_{2,t}
\stackrel{\mathrm{u.c.p.}} {\longrightarrow}0.
\]
For any fixed $x,y \in\mathbb{R}$ it holds that
\[
\widetilde{F}_n(t,x,y) \stackrel{\mathrm{u.c.p.}} {\longrightarrow}
\widetilde{F}(t,x,y):=\int_0^t \Phi
_{\sigma_s}(x)\Phi_{\sigma_s}(y) \,ds.
\]
This can be shown similarly to the proof of Proposition~\ref{prop1} as
follows. Let $\xi_j^n = n^{-1}\mathbh{1}_{\{\alpha_j^n \le x_1\}
}\mathbh{1}_{\{\alpha_{j+1}^{\prime n} \le x_2\}}$. Then
\[
\sum_{j=1}^{[nt]-1} \mathbb{E}\bigl[
\xi_j^n| \mathcal{F}_{(j-1)/{n}}\bigr] =
\frac{1}{n}\sum_{j=1}^{[nt]-1}
\Phi_{\sigma_{(j-1)/{n}}}(x_1)\Phi_{\sigma_{(j-1)/{n}}}(x_2)
\stackrel{\mathrm{u.c.p.}} {\longrightarrow}\widetilde{F}(t,x,y).
\]
On the other hand, we trivially have that $\sum_{j=1}^{[nt]-1} \mathbb
{E}[|\xi_j^n|^2| \mathcal{F}_{(j-1)/{n}}] \stackrel{\mathbb
{P}}{\longrightarrow}0$,
for any fixed $t>0$. Hence, the Lenglart's domination property (see
\cite{JS}, page~35) implies the convergence
\[
\sum_{j=1}^{[nt]-1} \bigl(
\xi_j^n - \mathbb{E}\bigl[\xi_j^n|
\mathcal{F}_{(j-1)/{n}}\bigr] \bigr) \stackrel{\mathrm{u.c.p.}} {
\longrightarrow}0,
\]
which in turn means that $\widetilde{F}_n(t,x,y) \stackrel{\mathrm
{u.c.p.}}{\longrightarrow}\widetilde{F}(t,x,y)$.

We know now that $V_{2,t}^{\prime n}$ converges to the claimed limit if
$G_2$ is compactly supported. For a general $G_2$ with polynomial
growth one can proceed exactly as in Proposition~\ref{prop1}.
To complete the proof, one has to show that $V_{2,t}^n-V_{2,t}^{\prime
n} \stackrel{\mathrm{u.c.p.}}{\longrightarrow}0$. This works exactly
as in Section~\ref{subsecsub31}.
\end{pf}
The properties of stable convergence immediately imply the following theorem.
%
%
\begin{theorem} \label{feasclt}
Let the assumptions of Theorem~\ref{th2} be satisfied. Let $t > 0$ be
fixed. Then we obtain the standard central limit theorem
%
%
\begin{equation}
\label{eqstandardclt} \frac{\sqrt{n} (U(H)_t^n - U(H)_t )}{\sqrt
{V_t^n}} \stackrel{d} {\longrightarrow}\mathcal{N}(0,1),
\end{equation}
where $V_t^n = V_{1,t}^n - V_{2,t}^n$ using the notation defined above.
\end{theorem}
The convergence in law in (\ref{eqstandardclt}) is a feasible central
limit theorem that can be used in statistical applications. It is
possible to obtain similar multivariate central limit theorems for
finite-dimensional vectors $\sqrt{n} (U(H)_{t_j}^n -
U(H)_{t_j} )_{1\le j \le k}$; we leave the details to the
interested reader.

\section{Statistical applications}\label{secstatex}
In this section we present some statistical applications of the limit
theory for $U$-statistics of continuous It\^o semimartingales.

\subsection{Gini's mean difference}

Gini's mean difference is a classical measure of statistical dispersion,
which serves as robust measure of variability of a probability
distribution \cite{David1968}. Recall that for a given distribution
$\mathbb Q$, Gini's mean difference is defined as
\[
\MD:=\mathbb{E}\bigl[|Y_1-Y_2|\bigr],
\]
where $Y_1,Y_2$ are independent random variables with distribution
$\mathbb Q$. In the framework of i.i.d. observations
$(Y_i)_{i\geq1}$, the measure $\MD$ is consistently estimated by the
$U$-statistic $\frac{2}{n(n-1)} \sum_{1\leq i<j\leq n} |Y_i-Y_j|$.
Gini's\vspace*{1pt} mean difference is connected to questions of stochasic dominance
as shown by \cite{Yitzhaki1982}.
We refer to the recent paper~\cite{LL} for the estimation theory for
Gini's mean difference under long range dependence.

In the setting of continuous It\^o semimartingales we conclude by
Theorem~\ref{th1} that
\[
U(H)_t^n \stackrel{\mathrm{u.c.p.}} {
\longrightarrow}\MD_t:= m_1 \int_{[0,t]^2}
\bigl|\sigma_{s_1}^2 + \sigma_{s_2}^2\bigr|^{1/2}
\,ds_1 \,ds_2,
\]
where the function $H$ is given by $H(x,y)=|x-y|$, and $m_p$ is the
$p$th absolute moment of $\mathcal N(0,1)$.
In mathematical finance the quantity $\MD_t$ may be viewed as an
alternative measure of price variability, which
is more robust to outliers than the standard quadratic variation $[X,X]_t$.

Formally, we cannot directly apply Theorem~\ref{th2} to obtain a weak
limit theory for the statistic $U(H)_t^n$,
since the function $H(x,y)=|x-y|$ is not differentiable, and $H$ is not
even in \textit{each} component. Since $Y_1-Y_2$
and $Y_1+Y_2$ have the same distribution for centered independent
normally distributed random variables $Y_1,Y_2$, the modification
\[
\widebar H(x,y):=\tfrac{1}2 \bigl(|x-y|+|x+y|\bigr),
\]
which is even in each component, has the same limit, that is,
$U(\widebar H)_t^n \stackrel{\mathrm{u.c.p.}}{\longrightarrow
}\MD_t$. Moreover,
using sub-differential calculus and defining
\[
\operatorname{grad} \widebar H(x,y):=\tfrac{1}2 \bigl(\operatorname
{sign}(x-y)+\operatorname{sign}(x+y), \operatorname{sign}(x-y)+
\operatorname{sign}(x+y) \bigr),
\]
all the proof steps remain valid (we also refer to \cite{BGJPS}, who
prove the central limit theorem for nondifferentiable
functions). Thus, by the assertion of Theorem~\ref{th2}, we deduce the
stable convergence
\[
\sqrt{n} \bigl( U(\widebar H)_t^n - \MD_t
\bigr) \stackrel{\mathrm{st}} {\longrightarrow}L_t= \int_{\mathbb{R}^2}
\bigl(|x_1-x_2| + |x_1+x_2|\bigr) \mathbb
G(t,dx_1) F(t,dx_2),
\]
where the stochastic fields $\mathbb G(t,x)$ and $F(t,x)$ are defined
in Proposition~\ref{prop2} and~(\ref{fnconv}), respectively.
Now, we follow the route proposed in Section~\ref{seccondvar} to
obtain a standard central limit theorem. We compute
the symmetrization $\widetilde G_1, \widetilde G_2$ of the functions
$G_1, G_2$ defined at (\ref{eqG1}) and (\ref{eqG2}), respectively:
\begin{eqnarray*}
\widetilde G_1 (x_1,x_2,x_3)&=&
\tfrac{1}6 \bigl(\bigl(|x_1-x_2| +|x_1+x_2|\bigr)
\bigl(|x_1-x_3| +|x_1+x_3|\bigr)
\\
&&\hspace*{9pt}{}+ \bigl(|x_2-x_1| +|x_2+x_1|\bigr)
\bigl(|x_2-x_3| +|x_2+x_3|\bigr)
\\
&&\hspace*{9pt}{}+ \bigl(|x_3-x_1| +|x_3+x_1|\bigr)
\bigl(|x_3-x_2| +|x_3+x_2|\bigr) \bigr),
\\
\widetilde G_1 (x_1,x_2;y_1,y_2)&=&
\tfrac{1}4 \bigl( \bigl(|x_1-y_1|
+|x_1+y_1|\bigr) \bigl(|x_2-y_2|
+|x_2+y_2|\bigr)
\\
&&\hspace*{9pt}{}+ \bigl(|x_1-y_2| +|x_1+y_2|\bigr)
\bigl(|x_2-y_1| +|x_2+y_1|\bigr) \bigr).
\end{eqnarray*}
Using these functions we construct the statistics $V_{1,t}^n$ and
$V_{2,t}^n$ (see Section~\ref{seccondvar}). Finally, for any fixed
$t>0$ we obtain a feasible central limit theorem
\[
\frac{\sqrt{n} (U(\widebar H)_t^n - \MD_t )}{\sqrt
{V_{1,t}^n-V_{2,t}^n}} \stackrel{d} {\longrightarrow}\mathcal{N}(0,1).
\]
The latter enables us to construct confidence regions for mean
difference statistic~$\MD_t$.

\subsection{$\mathbb{L}^p$-type tests for constant volatility}
In this subsection we propose a new homoscedasticity test for the
volatility process $\sigma^2$.
Our main idea relies on a certain distance measure, which is related to
$\mathbb L^p$-norms; we refer
to \cite{DP,DPV} for similar testing procedures in the $\mathbb L^2$ case.
Let us define
\[
h(s_1,\ldots, s_d):= \sum_{i=1}^d
\sigma_{s_i}^2, \qquad s_1,\ldots,
s_d\in[0,1],
\]
and consider a real number $p>1$. Our test relies on the $\mathbb L^p$-norms
\[
\|h\|_{\mathbb L^p}:= \biggl( \int_{[0,1]^d} \bigl|h(
\mathbf{s})\bigr|^p \,d\mathbf{s} \biggr)^{1/p}.
\]
Observe the inequality $\|h\|_{\mathbb L^p}\geq\|h\|_{\mathbb L^1}$
and, when the process $h$ is continuous, equality
holds if and only if $h$ is constant. Applying this intuition, we
introduce a distance measure $\mathcal M^2$ via
\[
\mathcal M^2:= \frac{\|h\|_{\mathbb L^p}^p - \|h\|_{\mathbb L^1}^p}{\|
h\|_{\mathbb L^p}^p}\in[0,1].
\]
Notice that a continuous process $\sigma^2$ is constant if and only if
$\mathcal M^2=0$. Furthermore, the measure $\mathcal M^2$ provides a
quantitative account of the deviation from the homoscedasticity
hypothesis, as it takes values in $[0,1]$.

For simplicity of exposition
we introduce an empirical analogue of $\mathcal M^2$ in the case $d=2$.
We define the functions
\[
H_1(x):= \tfrac{1}2 \bigl(|x_1-x_2|^{2p}
+ |x_1+x_2|^{2p}\bigr), \qquad
H_2(x):= x_1^2+x_2^2
\]
with $x\in\mathbb{R}^2$. Notice that both functions are continuously
differentiable and even in each component;
hence they satisfy the\vadjust{\goodbreak} assumptions of Theorems~\ref{th1} and~\ref
{th2}. In particular, Theorem~\ref{th1}
implies the convergence in probability
\[
U(H_1)_1^n \stackrel{\mathbb{P}} {
\longrightarrow}U(H_1)_1= m_{2p} \| h
\|_{\mathbb L^p}^p, \qquad U(H_2)_1^n
\stackrel{\mathbb{P}} {\longrightarrow}U(H_2)_1= \|h
\|_{\mathbb L^1},
\]
where the constant $m_{2p}$ has been defined in the previous
subsection. The main ingredient for a formal testing procedure
is the following result.

%
\begin{proposition} \label{multivariate}
Assume that conditions of Theorem~\ref{th2} hold. Then we obtain the
stable convergence
%
%
\begin{eqnarray}
\label{clthom} && \sqrt n \bigl( U(H_1)_1^n
- m_{2p} \|h\|_{\mathbb
L^p}^p, U(H_2)_1^n
- \|h\|_{\mathbb L^1} \bigr)\nonumber
\\
&&\qquad \stackrel{\mathrm{st}} {\longrightarrow}2 \biggl(\int_{\mathbb{R}^2}
H_1(x_1,x_2) \mathbb G(1,dx_1)
F(1,dx_2),
\\
&&\hspace*{59pt} \int_{\mathbb{R}^2} H_2(x_1,x_2)
\mathbb G(1,dx_1) F(1,dx_2) \biggr).
\nonumber
\end{eqnarray}
Furthermore, the $\mathcal F$-conditional covariance matrix
$V=(V_{ij})_{1\leq i,j\leq2}$ of the limiting random variable
is given as
%
%
\begin{eqnarray}
\label{vdef} V_{ij}&=& \int_0^1
\biggl( \int_{\mathbb{R}} f_{i}(x) f_{j}(x)
\varphi_{\sigma_s} (x) \,dx
\nonumber\\[-8pt]\\[-8pt]
&&\hspace*{19pt}{}- \biggl(\int_{\mathbb{R}} f_{i}(x)
\varphi_{\sigma_s} (x) \,dx \biggr) \biggl(\int_{\mathbb{R}}
f_{j}(x) \varphi_{\sigma_s} (x) \,dx \biggr) \biggr) \,ds
\nonumber
\end{eqnarray}
with
\[
f_i(x):=2\int_{\mathbb{R}} H_i(x,y)
F(1,dy), \qquad i=1,2.
\]
\end{proposition}
\begin{pf}
As in the proof of Theorem~\ref{th2} we deduce that
\[
\sqrt{n} \bigl(U(H_i)_1^n -
U(H_i)_1\bigr) = L_1^{\prime n}(i) +
o_{\mathbb P} (1), \qquad i=1,2,
\]
where $L_1^{\prime n}(i)$ is defined via
\[
L_1^{\prime n}(i)= 2 \int_{\mathbb{R}^2}
H_i(x_1,x_2) \mathbb G_n (1,
dx_1) F(1, dx_2).
\]
Now, exactly as in steps (ii)--(v) of the proof of Proposition~\ref
{prop3}, we conclude the joint
stable convergence in (\ref{clthom}). The $\mathcal F$-conditional
covariance matrix $V$ is obtained
from Remark~\ref{rem2}(iii) as in the beginning of Section~\ref{seccondvar}.
\end{pf}
Let now $\mathcal M_n^2$ be the empirical analogue of $\mathcal M^2$,
that is,
\[
\mathcal M_n^2:= \frac{m_{2p}^{-1} U(H_1)_1^n -
(U(H_2)_1^n)^p}{m_{2p}^{-1} U(H_1)_1^n} \stackrel{\mathbb
{P}} {\longrightarrow}\mathcal M^2.
\]
Observe the identities
\[
\mathcal M_n^2 = r \bigl(U(H_1)_1^n,
U(H_2)_1^n \bigr), \qquad\mathcal
M^2= r \bigl(m_{2p}\|h\|_{\mathbb L^p}^p,
\|h\|_{\mathbb L^1} \bigr),\vadjust{\goodbreak}
\]
where $r(x,y)=1-m_{2p} \frac{y^p}{x}$. Applying Proposition~\ref
{multivariate} and delta method for stable convergence, we conclude
that $\sqrt{n}(\mathcal M_n^2 - \mathcal M^2)$ converges stably in law
toward a mixed normal distribution with mean $0$
and $\mathcal F$-conditional variance given by
\[
v^2:=\nabla r \bigl(m_{2p}\|h\|_{\mathbb L^p}^p,
\|h\|_{\mathbb L^1} \bigr) V \nabla r \bigl(m_{2p}\|h
\|_{\mathbb L^p}^p, \|h\|_{\mathbb L^1} \bigr)^{\star},
\]
where the random variable $V\in\mathbb R^{2 \times2}$ is defined at
(\ref{vdef}).

For an estimation of $V$ we can proceed as in Section~\ref
{seccondvar}. Define the functions $G_1^{ij}\dvtx\mathbb{R}^3 \to
\mathbb
{R}$ and $G_2^{ij}\dvtx \mathbb{R}^4 \to\mathbb{R}$ by
\begin{eqnarray*}
G_1^{ij}(x_1,x_2,x_3)&=&H_i(x_1,x_2)H_j(x_1,x_3),
\\
G_2^{ij}(x_1,x_2,y_1,y_2)&=&H_i(x_1,y_1)H_j(x_2,y_2),
\qquad i,j=1,2.
\end{eqnarray*}
Let further $\widetilde{G}_1^{ij}$ be the symmetrization of $G_1^{ij}$
and $\widetilde{G}_2^{ij}$ the symmetrization of $G_2^{ij}$ with
respect to the $\mathbf{y}$-values. With
\[
W_{ij}:=\frac{4}{n}\pmatrix{n
\cr
2}^{-1} \sum
_{i_1=1}^{n-1} \sum_{1
\leq i_2 < i_3 \leq n}
\widetilde{G}_2^{ij}\bigl(\sqrt{n} \Delta_{i_1}^n
X,\sqrt{n} \Delta_{i_1+1}^n X,\sqrt{n} \Delta_{i_2}^n
X,\sqrt{n} \Delta_{i_3}^n X\bigr),
\]
we can, exactly as in Section~\ref{seccondvar}, deduce that
\[
V^n:= \bigl(4U\bigl(\widetilde{G}_1^{ij}
\bigr)_1^n-W_{ij}\bigr)_{i,j=1,2}
\stackrel{\mathbb{P}} {\longrightarrow}V.
\]
Using the previous results, we directly get
\[
v_n^2:= \nabla r \bigl(U(H_1)_1^n,
U(H_2)_1^n \bigr) V^n \nabla r
\bigl(U(H_1)_1^n, U(H_2)_1^n
\bigr)^{\star} \stackrel{\mathbb{P}} {\longrightarrow}v^2.
\]
Now the properties of stable convergence yield the following feasible
central limit theorem:
%
%
\begin{equation}
\label{eqlpclt} \frac{\sqrt{n}(\mathcal M_n^2 - \mathcal M^2)}{\sqrt
{v_n^2}} \stackrel{d} {\longrightarrow}\mathcal N(0,1).
\end{equation}
With these formulas at hand, we can derive a formal test procedure for
the hypothesis
\[
H_0\dvtx\sigma_s^2 \mbox{ is constant on }
[0,1] \quad\mbox{vs.}\quad H_1\dvtx\sigma_s^2\mbox{ is not constant
on } [0,1].
\]
These hypotheses are obviously equivalent to
\[
H_0\dvtx\mathcal M^2=0\quad\mbox{vs.}\quad
H_1\dvtx\mathcal M^2>0.
\]
Defining the test statistic $S_n$ via
\[
S_n:= \frac{\sqrt{n} \mathcal M_n^2}{\sqrt{v_n^2}},
\]
we reject the null hypothesis at level $\gamma\in(0,1)$ whenever
$S_t^n>c_{1-\gamma}$, where
$c_{1-\gamma}$ denotes the $(1-\gamma)$-quantile of $\mathcal
N(0,1)$. Now, (\ref{eqlpclt})
implies that
\[
\lim_{n\rightarrow\infty} \mathbb P_{H_0} (S_n>c_{1-\gamma})
= \gamma, \qquad\lim_{n\rightarrow\infty} \mathbb P_{H_1}
\bigl(S_n^n>c_{1-\gamma}\bigr) = 1.\vadjust{\goodbreak}
\]
In other words, our test statistic is consistent and keeps the level
$\gamma$ asymptotically.

\subsection{Wilcoxon test statistic for structural breaks}
Change-point analysis has been an active area of research for many
decades; we refer to \cite{CH} for a comprehensive
overview.
The Wilcoxon statistic is a standard statistical procedure for testing
structural breaks in location models.
Let $(Y_i)_{1\leq i\leq n}$, $(Z_i)_{1\leq i\leq m}$ be mutually
independent observations with
$Y_i\sim\mathbb Q_{\theta_1}$, $Z_i\sim\mathbb Q_{\theta_2}$, where
$\mathbb Q_{\theta}(A)= \mathbb Q_{0}(A-\theta)$
for all $A\in\mathcal B(\mathbb{R})$ and $\mathbb Q_{0}$ be a
nonatomic probability measure. In this classical framework
the Wilcoxon statistic is defined~by
\[
\frac{1}{nm} \sum_{i=1}^{n} \sum
_{j=1}^m \mathbh{1}_{\{Y_i\leq Z_j\}}.
\]
Under the null hypothesis $\theta_1=\theta_2$, the test statistic is
close to $1/2$, while deviations from this value
indicate that $\theta_1 \neq\theta_2$. We refer to the recent work
\cite{DRT} for change-point tests for long-range
dependent data.

Applying the same intuition we may provide a test statistic for
structural breaks in the volatility process $\sigma^2$.
Assume that the semimartingale $X$ is observed at high frequency on the
interval $[0,1]$ and the volatility
is constant on the intervals $[0,t)$ and $(t,1]$ for some $t\in(0,1)$,
that is, $\sigma_s^2=\sigma_0^2$ on $[0,t)$
and $\sigma_s^2=\sigma_1^2$ on $(t,1]$. Our aim is to test the null
hypothesis $\sigma_0^2=\sigma_1^2$
or to infer the change-point $t$ when $\sigma_0^2 \neq\sigma_1^2$.
In this framework
the Wilcoxon type statistic is defined via
\[
\WL_t^n:= \frac{1}{n^2}\sum
_{i=1}^{[nt]} \sum_{j=[nt]+1}^n
\mathbh{1}_{\{|\Delta_i^n X|\leq|\Delta_j^n X| \}}.
\]
Notice that the kernel is neither symmetric nor continuous.
Nevertheless, we deduce the following result.

%
\begin{proposition} \label{prop}
Assume that condition (\ref{nonvanish}) holds. Then we obtain the convergence:
%
%
\begin{eqnarray}
\label{wil} \WL_t^n \stackrel{\mathrm{u.c.p.}} {
\longrightarrow}\WL_t&:=& \int_{0}^{t}\!\!\int_{t}^1 \biggl( \int_{\mathbb
{R}^2}
\mathbh{1}_{\{|\sigma_{s_1} u_1|\leq|\sigma_{s_2} u_2| \}} \varphi
_d(\mathbf{u}) \,d\mathbf{u} \biggr)
\,ds_1\,ds_2
\\
&=&\int_{0}^{t}\!\! \int_{t}^1
\biggl( 1-\frac{2}{\pi}\arctan\biggl\llvert\frac{\sigma_{s_1}}{\sigma
_{s_2}}\biggr\rrvert
\biggr) \,ds_1\,ds_2.
\end{eqnarray}
\end{proposition}

\begin{pf}
As in the proof of Theorem~\ref{th1}, we first show
the convergence (\ref{wil}) for the approximations $\alpha_i^n$ of
the scaled increments $\sqrt n \Delta_i^n X$. We define
\[
U_t^{\prime n}:= \int_{\mathbb{R}^2}
\mathbh{1}_{\{ |x|\leq|y|\}} F_n(t,dx) \bigl(F_n(1,dy)-F_n(t,dy)
\bigr).\vadjust{\goodbreak}
\]
Since condition (\ref{nonvanish}) holds, the measure $F_n(t,dx)$ is
nonatomic. Hence, we conclude that
\[
U_t^{\prime n} \stackrel{\mathrm{u.c.p.}} {\longrightarrow}\WL_t
\]
exactly as in the proof of Proposition~\ref{prop1}. It remains to
prove the convergence
\[
\WL_t^n - U_t^{\prime n} \stackrel{
\mathrm{u.c.p.}} {\longrightarrow}0.
\]
Observe the identity
\begin{eqnarray*}
\WL_t^n - U_t^{\prime n} &=&
\frac{1}{n^2}\sum_{i=1}^{[nt]} \sum
_{j=[nt]+1}^n ( \mathbh{1}_{\{|\Delta_i^n X|\leq|\Delta_j^n X| \}} -
\mathbh{1}_{\{|\alpha_i^n|\leq|\alpha_j^n| \}} )
\\
&=& \frac{1}{n^2}\sum_{i=1}^{[nt]} \sum
_{j=[nt]+1}^n ( \mathbh{1}_{\{|\Delta_i^n X|\leq|\Delta_j^n X| \}} -
\mathbh{1}_{\{|\sqrt{n} \Delta_i^n X|\leq|\alpha_j^n| \}}
\\
&&\hspace*{73pt}{}+ \mathbh{1}_{\{|\sqrt{n} \Delta_i^n
X|\leq|\alpha_j^n| \}} -
\mathbh{1}_{\{|\alpha_i^n|\leq|\alpha_j^n| \}} ).
\end{eqnarray*}
In the following we concentrate on proving that
\[
\frac{1}{n^2}\sum_{i=1}^{[nt]} \sum
_{j=[nt]+1}^n ( \mathbh{1}_{\{|\sqrt{n} \Delta_i^n X|\leq|\alpha_j^n| \}
} -
\mathbh{1}_{\{|\alpha_i^n|\leq|\alpha_j^n| \}} ) \stackrel{\mathrm
{u.c.p.}} {\longrightarrow}0,
\]
as the other part is negligible by the same arguments. Using the identity
\begin{eqnarray*}
&& \mathbh{1}_{\{ |\sqrt{n} \Delta_i^n X|\leq|\alpha_j^n| \}} -\mathbh
{1}_{\{|\alpha_i^n|\leq|\alpha_j^n| \}}
\\
&&\qquad = \mathbh{1}_{\{
|\sqrt{n} \Delta_i^n X|\leq|\alpha_j^n|, |\alpha_i^n|> |\alpha
_j^n| \}} -\mathbh{1}_{\{|\sqrt{n} \Delta_i^n X|> |\alpha_j^n|, |\alpha
_i^n|\leq|\alpha_j^n| \}}
\end{eqnarray*}
we restrict our attention on proving
\[
\frac{1}{n^2}\sum_{i=1}^{[nt]} \sum
_{j=[nt]+1}^n \mathbh{1}_{\{
|\sqrt{n} \Delta_i^n X|> |\alpha_j^n|, |\alpha_i^n|\leq|\alpha
_j^n| \}}
\stackrel{\mathrm{u.c.p.}} {\longrightarrow}0.
\]
For an arbitrary $q\in(0,1/2)$, we deduce the inequality
\begin{eqnarray*}
\mathbb{E}[\mathbh{1}_{\{|\sqrt{n} \Delta_i^n X|> |\alpha_j^n|,
|\alpha_i^n|\leq|\alpha_j^n| \}}] &\leq&\mathbb{E} \biggl[
\frac{|\sqrt{n} \Delta_i^n X-\alpha
_i^n|^q}{||\alpha_j^n| - |\alpha_i^n||^q} \biggr]
\\
&\leq&\mathbb{E}\bigl[\bigl|\sqrt{n} \Delta_i^n X-
\alpha_i^n\bigr|^{2q}\bigr]^{1/2}
\mathbb{E}\bigl[\bigl|\bigl|\alpha_j^n\bigr| - \bigl|\alpha_i^n\bigr|\bigr|^{-2q}
\bigr]^{1/2}.
\end{eqnarray*}
For a standard normal random variable $U$, and for any $x>0, y\geq0$, define
\[
g_q(x,y):=\mathbb{E}\bigl[\bigl|x|U|-y\bigr|^{-2q}\bigr].
\]
Since $2q <1$, we have
%
%
\begin{eqnarray}\label{ndin}
\qquad g_q(x,y) &=& \mathbb{E}\bigl[\bigl|x|U|-y\bigr|^{-2q}
\mathbh{1}_{\{|x|U|-y|\leq1 \}}\bigr]
+\mathbb{E}\bigl[\bigl|x|U|-y\bigr|^{-2q}
\mathbh{1}_{\{|x|U|-y|> 1 \}}\bigr]\nonumber
\\
&\leq&\int_{\mathbb{R}}\bigl|x|u|-y\bigr|^{-2q}\mathbh{1}_{\{|x|u|-y|\leq1 \}
}\,du +1
\\
&\leq&\frac{C_q}{x}+1 < \infty.\nonumber
\end{eqnarray}
Due to assumption (\ref{nonvanish}) and by a localization argument, we
can assume that $\sigma_t$ is uniformly bounded away from zero.
Therefore, and
by (\ref{ndin}) we obtain
\begin{eqnarray*}
\mathbb{E}\bigl[\bigl|\bigl|\alpha_j^n\bigr| - \bigl|\alpha_i^n\bigr|\bigr|^{-2q}
\bigr] &=& \mathbb{E}\bigl[\mathbb{E}\bigl[ \bigl|\bigl|\alpha_j^n\bigr|
- \bigl|\alpha_i^n\bigr|\bigr|^{-2q}| \mathcal
{F}_{(j-1)/{n}}\bigr]\bigr]
\\
&=& \mathbb{E}\bigl[ g_q\bigl(\sigma_{(j-1)/{n}},
\alpha_i^n\bigr)\bigr] \leq C_q < \infty.
\end{eqnarray*}
Hence
\begin{eqnarray*}
&&\frac{1}{n^2}\sum_{i=1}^{[nt]} \sum
_{j=[nt]+1}^n \mathbb{E}[\mathbh{1}_{\{\bigl|\sqrt{n} \Delta_i^n X\bigr|> |\alpha
_j^n|, |\alpha
_i^n|\leq|\alpha_j^n| \}}]
\\
&&\qquad \leq\frac{C}{n^2}\sum_{i=1}^{[nt]}
\sum_{j=[nt]+1}^n \mathbb{E}\bigl[\bigl|\sqrt{n}
\Delta_i^n X-\alpha_i^n\bigr|^{2q}
\bigr]^{1/2} \stackrel{\mathrm{u.c.p.}} {\longrightarrow}0,
\end{eqnarray*}
where the last convergence follows as in (\ref{plln}). This completes
the proof of Proposition~\ref{prop}.
\end{pf}

Now, observe that when the process $\sigma^2$ has no change-point at
time $t\in(0,1)$ (i.e., $\sigma_0^2=\sigma_1^2$)
the limit at (\ref{wil}) is given by $\WL_t=\frac{1}2t(1-t)$. Thus,
under the null hypothesis $\sigma_0^2=\sigma_1^2$,
we conclude that $\WL_t^n \stackrel{\mathrm{u.c.p.}}{\longrightarrow
}\frac{1}2 t(1-t)$. Since the time point
$t\in(0,1)$ is unknown in general, we may use the test
statistic
\[
\sup_{t\in(0,1)} \biggl|\WL_t^n -
\frac{1}2t(1-t)\biggr|
\]
to test for a possible change point. Large values of this quantity
speak against the null hypothesis.
On the other hand, under the alternative $\sigma_0^2 \neq\sigma
_1^2$, the statistic $\hat{t}_n:=
\operatorname{argsup}_{t\in(0,1)} |\WL_t^n - \frac{1}2t(1-t)|$ provides
a consistent estimator of the change-point $t\in(0,1)$.
A formal testing
procedure would rely on a stable central limit theorem for $\WL_t^n$,
which is expected to be highly complex,
since the applied kernel is not differentiable.

\section{Proofs of some technical results} \label{proofs}

Before we start with the proofs of (\ref{sec3a}) and (\ref
{sec3b}) we state the following lemma, which can be shown exactly as
\cite{BGJPS}, Lemma~5.4.\vadjust{\goodbreak}
%
%
\begin{lemma} \label{lem2}
Let $f\dvtx\mathbb{R}^{d} \to\mathbb{R}^{q}$ be a continuous function of
polynomial growth. Let further $\gamma_i^n, \gamma_i'^n$ be
real-valued random variables satisfying $\mathbb{E}[(|\gamma_i^n|+
|\gamma_i'^n|)^p]$ $\leq C_p$ for all $p\geq2$ and
\[
\pmatrix{n
\cr
d}^{-1} \sum_{\mathbf{i} \in\mathcal{A}_t^n(d)}
\mathbb{E}\bigl[\bigl\|\gamma_{\mathbf{i}}^n- \gamma_{\mathbf{i}}'^n
\bigr\|^2\bigr] \to0.
\]
Then we have for all $t>0$,
\[
\pmatrix{n
\cr
d}^{-1} \sum_{\mathbf{i} \in\mathcal{A}_t^n(d)}
\mathbb{E}\bigl[\bigl\|f\bigl(\gamma_{\mathbf{i}}^n\bigr)-f\bigl(
\gamma_{\mathbf{i}}'^n\bigr)\bigr\| ^2\bigr]
\to0.
\]
\end{lemma}

Recall that we assume (\ref{probound1}) without loss of generality; in
Sections~\ref{subsecsub32}~and~\ref{subsecsub33} we further
assume (\ref{probound2}), that is, all the involved processes are bounded.

\subsection{Proof of (\texorpdfstring{\protect\ref{sec3a}}{13})}\label{subsecsub31}
The Burkh\"older inequality yields that $\mathbb{E}[(|\sqrt{n}\Delta
_i^n X|+ |\alpha_i^n|)^p] \leq C_p$ for all $p\geq2$. In view of the
previous lemma $U(H)^n-\widetilde{U}(H)^n \stackrel{\mathrm
{u.c.p.}}{\longrightarrow}0$ is a direct
consequence of
%
%
\begin{equation}
\label{plln} \qquad \pmatrix{n
\cr
d}^{-1} \sum
_{\mathbf{i} \in\mathcal{A}_t^n(d)} \mathbb{E}\bigl[\bigl\|\sqrt{n}\Delta
_{\mathbf{i}}^n
X - \alpha_{\mathbf
{i}}^n\bigr\|^2\bigr] \leq
\frac{C}{n}\sum_{j=1}^{[nt]} \mathbb{E}
\bigl[\bigl|\sqrt{n}\Delta_{j}^n X - \alpha_{j}^n\bigr|^2
\bigr] \to0
\end{equation}
as it is shown in \cite{BGJPS}, Lemma 5.3.

\subsection{Proof of (\texorpdfstring{\protect\ref{sec3b}}{35})}\label{subsecsub32}
We divide the proof into several steps.
\begin{longlist}[(iii)]
\item[(i)] We claim that
\[
\sqrt{n}\bigl(U(H)^n-\widetilde{U}(H)^n
\bigr)-P^n(H) \stackrel{\mathrm{u.c.p.}} {\longrightarrow}0,
\]
where
\[
P_t^n(H):=\sqrt{n}\pmatrix{n
\cr
d}^{-1} \sum
_{\mathbf{i} \in
\mathcal{A}_t^n(d)} \nabla H\bigl(\alpha_{\mathbf{i}}^n
\bigr) \bigl(\sqrt{n}\Delta_{\mathbf{i}}^n X -
\alpha_{\mathbf{i}}^n\bigr).
\]
Here, $\nabla H$ denotes the gradient of $H$. This can be seen as follows.
Since the process $\sigma$ is itself a continuous It\^o
semimartingale, we have
%
%
\begin{equation}
\label{est1} \mathbb{E}\bigl[\bigl|\sqrt{n}\Delta_i^n X-
\alpha_i^n\bigr|^p\bigr] \leq C_p
n^{-p/2}
\end{equation}
for all $p\geq2$.
By the mean value theorem, for any $\mathbf{i}\in\mathcal
{A}_t^n(d)$, there exists a random variable $\chi_{\mathbf{i}}^n \in
\mathbb{R}^d$ such that
\[
H\bigl(\sqrt{n}\Delta_{\mathbf{i}}^n X\bigr) - H\bigl(
\alpha_{\mathbf
{i}}^n\bigr)=\nabla H\bigl(\chi_{\mathbf{i}}^n
\bigr) \bigl(\sqrt{n}\Delta_{\mathbf
{i}}^n X-\alpha_{\mathbf{i}}^n
\bigr)
\]
with $\| \chi_{\mathbf{i}}^n-\alpha_{\mathbf{i}}^n\|\leq\| \sqrt
{n}\Delta_{\mathbf{i}}^n X -\alpha_{\mathbf{i}}^n\|$. Therefore, we have
\begin{eqnarray*}
&& \mathbb{E}\Bigl[\sup_{t\leq T} \bigl|\sqrt{n}\bigl(U(H)_t^n-
\widetilde{U}_t(H)^n\bigr)-P_t^n(H)\bigr|
\Bigr]
\\
&&\qquad \leq C\sqrt{n}\pmatrix{n
\cr
d}^{-1} \sum
_{\mathbf{i} \in\mathcal
{A}_T^n(d)} \mathbb{E}\bigl[\bigl\|(\nabla H\bigl(\chi_{\mathbf{i}}^n
\bigr)-\nabla H\bigl(\alpha_{\mathbf{i}}^n\bigr)\bigr\| \bigl\|\bigl(\sqrt{n}
\Delta_{\mathbf{i}}^n X - \alpha_{\mathbf{i}}^n\bigr)
\bigr\|\bigr]
\\
&&\qquad \leq C\sqrt{n}\pmatrix{n
\cr
d}^{-1} \biggl( \sum
_{\mathbf{i} \in
\mathcal{A}_T^n(d)}\mathbb{E}\bigl[ \bigl\|\bigl(\nabla H\bigl(
\chi_{\mathbf
{i}}^n\bigr)-\nabla H\bigl(\alpha_{\mathbf{i}}^n
\bigr)\bigr)\bigr\|^2 \bigr] \biggr)^{1/2}
\\
&&\quad\qquad{} \times\biggl(\sum_{\mathbf{i} \in\mathcal
{A}_T^n(d)}\mathbb{E}\bigl[ \bigl\|
\bigl(\sqrt{n}\Delta_{\mathbf{i}}^n X - \alpha_{\mathbf{i}}^n
\bigr)\bigr\|^2\bigr] \biggr)^{1/2}
\\
&&\qquad \leq C \biggl\{\pmatrix{n
\cr
d}^{-1}\sum
_{\mathbf{i} \in\mathcal
{A}_T^n(d)}\mathbb{E}\bigl[ \bigl\|\bigl(\nabla H\bigl(
\chi_{\mathbf{i}}^n\bigr)-\nabla H\bigl(\alpha_{\mathbf{i}}^n
\bigr)\bigr)\bigr\|^2\bigr] \biggr\}^{1/2}
\\
&&\qquad \to0
\end{eqnarray*}
by (\ref{plln}) and Lemma~\ref{lem2}.

\item[(ii)] In this and the next step we assume that $H$ has compact support.
Now we split $P_t^n$ up into two parts:
%
%
\begin{eqnarray}
\label{spl1} P_t^n &=& \sqrt{n}\pmatrix{n
\cr
d}^{-1} \sum_{\mathbf{i} \in\mathcal
{A}_t^n(d)} \nabla H\bigl(
\alpha_{\mathbf{i}}^n\bigr)v_{\mathbf{i}}^n(1)
\nonumber\\[-8pt]\\[-8pt]
&&{} +
\sqrt{n}\pmatrix{n
\cr
d}^{-1} \sum_{\mathbf{i} \in\mathcal
{A}_t^n(d)}
\nabla H\bigl(\alpha_{\mathbf{i}}^n\bigr)v_{\mathbf{i}}^n(2),\nonumber
\end{eqnarray}
where $\sqrt{n}\Delta_{\mathbf{i}}^n X - \alpha_{\mathbf{i}}^n =
v_{\mathbf{i}}^n(1)+v_{\mathbf{i}}^n(2)$
and $\mathbf{i}=(i_1, \ldots, i_d)$, with
\begin{eqnarray*}
v_{i_k}^n(1) &=&\sqrt{n} \biggl(n^{-1}
a_{({i_k-1})/{n}}
+\int_{({i_k-1})/{n}}^{({i_k})/{n}} \bigl\{\tilde{\sigma
}_{({i_k-1})/{n}} (W_s-W_{({i_k-1})/{n}})
\\
&&\hspace*{131pt} {}+\tilde{v}_{({i_k-1})/{n}}(V_s-V_{({i_k-1})/{n}})
\bigr\} \,dW_s \biggr),
\\
v_{i_k}^n(2)&=&\sqrt{n} \biggl( \int_{({i_k-1})/{n}}^{({i_k})/{n}}(a_{s}-a_{({i_k-1})/{n}})\,ds
\\
&&\hspace*{20pt}{} +\int_{({i_k-1})/{n}}^{({i_k})/{n}} \biggl\{\int
_{({i_k-1})/{n}}^{s}\tilde{a}_u \,du
\\
&&\hspace*{75pt}{}
+ \int_{({i_k-1})/{n}}^{s}(\tilde{
\sigma}_{u-} -\tilde{\sigma}_{({i_k-1})/{n}}) \,dW_u
\\
&&\hspace*{79pt}{} +\int
_{({i_k-1})/{n}}^{s}(\tilde{v}_{u-}-
\tilde{v}_{({i_k-1})/{n}}) \,dV_u \biggr\} \,dW_s \biggr).
\end{eqnarray*}
We denote the first and the second summand on the right-hand side of
(\ref{spl1}) by $S_t^n$ and $\widetilde{S}_t^n$, respectively.
First, we show the convergence $\widetilde{S}{}^n \stackrel{\mathrm
{u.c.p.}}{\longrightarrow}0$. Since the
first derivative of $H$ is of polynomial growth we have $\mathbb{E}[\|
\nabla H (\alpha_{\mathbf{i}}^n)\|^2] \leq C$ for all $\mathbf{i}
\in\mathcal{A}_t^n(d)$.
Furthermore, we obtain by using the H\"older, Jensen and Burkh\"older
inequalities
\begin{eqnarray*}
&& \mathbb{E}\bigl[\bigl|v_{i_k}^n (2)\bigr|^2\bigr]
\\
&&\qquad \leq
\frac{C}{n^2}+\int_{({i_k-1})/{n}}^{({i_k})/{n}} (a_s
- a_{{[ns]}/{n}})^2+(\tilde{\sigma}_{s-} - \tilde{
\sigma}_{{[ns]}/{n}})^2+(\tilde{v}_{s-} -
\tilde{v}_{{[ns]}/{n}})^2 \,ds.
\end{eqnarray*}
Thus, for all $t>0$, we have
\begin{eqnarray*}
&&\sqrt{n}\pmatrix{n
\cr
d}^{-1}\mathbb{E}\sum
_{\mathbf{i} \in
\mathcal{A}_t^n(d)} \bigl| \nabla H\bigl(\alpha_{\mathbf{i}}^n
\bigr)v_{\mathbf
{i}}^n(2) \bigr|
\\
&&\qquad \leq C\sqrt{n}\pmatrix{n
\cr
d}^{-1} \biggl(\mathbb{E} \biggl[ \sum
_{\mathbf{i} \in\mathcal{A}_t^n(d)} \bigl\|\nabla H\bigl(\alpha_{\mathbf
{i}}^n
\bigr)\bigr\|^2 \biggr] \biggr)^{1/2} \biggl(\mathbb{E} \biggl[
\sum_{\mathbf
{i} \in\mathcal{A}_t^n(d)} \bigl\|v_{\mathbf{i}}^n(2)
\bigr\|^2 \biggr] \biggr)^{1/2}
\\
&&\qquad \leq C \Biggl( n \pmatrix{n
\cr
d}^{-1} \mathbb{E} \Biggl[ \sum
_{i_1,\ldots,i_d=1}^{[nt]} \bigl(\bigl|v_{i_1}^n(2)\bigr|^2+
\cdots+\bigl|v_{i_d}^n(2)\bigr|^2\bigr) \Biggr]
\Biggr)^{1/2}
\\
&&\qquad \leq C \Biggl(\mathbb{E} \Biggl[ \sum_{j=1}^{[nt]}
\bigl|v_{j}^n(2)\bigr|^2 \Biggr] \Biggr)^{1/2}
\\
&&\qquad \leq C \biggl(n^{-1}+\int_{0}^{t}
(a_s - a_{{[ns]}/{n}})^2+(\tilde{
\sigma}_{s-} - \tilde{\sigma}_{{[ns]}/{n}})^2+(
\tilde{v}_{s-} - \tilde{v}_{{[ns]}/{n}})^2 \,ds
\biggr)^{1/2}
\\
&&\qquad \to0
\end{eqnarray*}
by the dominated convergence theorem, and $\widetilde{S}{}^n \stackrel
{\mathrm{u.c.p.}}{\longrightarrow}0$
readily follows.

\item[(iii)] To show $S^n \stackrel{\mathrm{u.c.p.}}{\longrightarrow}0$ we use
\[
S_t^n = \sum_{k=1}^d
\sqrt{n}\pmatrix{n
\cr
d}^{-1}\sum_{\mathbf{i}
\in\mathcal{A}_t^n(d)}
\partial_k H\bigl(\alpha_{\mathbf
{i}}^n
\bigr)v_{i_k}^n(1)=:\sum_{k=1}^d
S_t^n(k).
\]
Before we proceed with proving $S^n(k) \stackrel{\mathrm
{u.c.p.}}{\longrightarrow}0$, for $k=1,\ldots,d$, we
make two observations: first, by the Burkh\"older inequality, we deduce
%
%
\begin{equation}
\label{v1} \mathbb{E}\bigl[\bigl|\sqrt{n} v_{i_k}^n(1)\bigr|^p
\bigr] \leq C_p\qquad\mbox{for all $p \ge2$,}
\end{equation}
and second, for fixed $x \in\mathbb{R}^{d-k}$, and for all $\mathbf
{i}=(i_1,\ldots,i_k)\in\mathcal{A}_t^n(k)$, we have
%
%
\begin{equation}
\label{cvan} \mathbb{E}\bigl[\partial_k H\bigl(
\alpha_{\mathbf{i}}^n, x\bigr)v_{i_k}^n(1) |
\mathcal{F}_{({i_k-1})/{n}}\bigr]=0,
\end{equation}
since $\partial_k H$ is an odd function in its $k$th component. Now,
we will prove that
%
%
\begin{equation}
\label{fcon} \sqrt{n}n^{-k}\sum_{\mathbf{i}\in\mathcal{A}_t^n(k)}
\partial_k H\bigl(\alpha_{\mathbf{i}}^n,x
\bigr)v_{i_k}^n(1) \stackrel{\mathrm{u.c.p.}} {
\longrightarrow}0,
\end{equation}
for any fixed $x \in\mathbb{R}^{d-k}$. From (\ref{cvan}) we know
that it suffices to show that
\[
\sum_{i_k=1}^{[nt]} \mathbb{E} \biggl[ \biggl(
\sum_{1\leq i_1<\cdots
<i_{k-1}<i_k} \chi_{i_1,\ldots, i_{k}}
\biggr)^2 \bigg| \mathcal{F}_{({i_k-1})/{n}} \biggr] \stackrel{\mathbb{P}}
{\longrightarrow}0,
\]
where $\chi_{i_1,\ldots, i_{k}}:= \sqrt{n}n^{-k} \partial_k
H(\alpha_{\mathbf{i}}^n,x)v_{i_k}^n(1)$. (Note that the sum in the
expectation only runs over the indices $i_1,\ldots,i_{k-1}$.) But this follows
from the $L^1$ convergence and (\ref{v1}) via
\begin{eqnarray*}
&& \sum_{i_k=1}^{[nt]} \mathbb{E} \biggl[
\biggl( \sum_{1\leq i_1<\cdots
<i_{k-1}<i_k} \chi_{i_1,\ldots, i_{k}}
\biggr)^2 \biggr]
\\
&&\qquad \leq\frac{C}{n^{k}} \sum_{i_k=1}^{[nt]}
\sum_{1\leq i_1<\cdots
<i_{k-1}<i_k} \mathbb{E} \bigl[\bigl(
\partial_k H\bigl(\alpha_{\mathbf
{i}}^n,x
\bigr)v_{i_k}^n(1) \bigr)^2 \bigr]
\\
&&\qquad \leq\frac{C}{n} \to0.
\end{eqnarray*}
Recall that we still assume that $H$ has compact support. Let
the support of $H$ be a subset of $[-K,K]^d$ and further $-K=z_0<\cdots
<z_m=K$ be an equidistant partition of $[-K,K]$. We denote the set
$ \{z_0,\ldots,z_m \}$ by $Z_m$. Also, let
$\eta(\varepsilon):=\sup\{\| \nabla H(\mathbf{x})- \nabla
H(\mathbf{y})\|; \|\mathbf{x}-\mathbf{y}\|\leq\varepsilon\}$
be the modulus of continuity of $\nabla H$. Then we have
\begin{eqnarray*}
\sup_{t \leq T} \bigl|S_t^n (k)\bigr| &\leq& C
\sqrt{n}n^{-k}\sup_{t \leq T} \sup_{x\in[-K,K]^{d-k}}
\Biggl| \sum_{\mathbf{i}\in\mathcal
{A}_t^n(k)} \partial_k H\bigl(
\alpha_{\mathbf{i}}^n,x\bigr)v_{i_k}^n(1) \Biggr|
\\
&\leq& C\sqrt{n} n^{-k} \sup_{t\leq T}\max
_{x \in Z_m^{d-k}} \Biggl| \sum_{\mathbf{i}\in\mathcal{A}_t^n(k)}
\partial_k H\bigl(\alpha_{\mathbf{i}}^n,x
\bigr)v_{i_k}^n(1) \Biggr|
\\
&&{} +C \sqrt{n} n^{-k}\sum_{\mathbf{i}\in\mathcal{A}_T^n(k)}\eta
\biggl(\frac{2K}{m} \biggr)\bigl|v_{i_k}^n(1)\bigr|.
\end{eqnarray*}
Observe that, for fixed $m$, the first summand converges in probability
to $0$ as \mbox{$n \to\infty$} by (\ref{fcon}). The second summand is
bounded in expectation by $C\eta(2K/m)$ which converges to $0$ as
$m\to\infty$. This implies $S_t^n(k) \stackrel{\mathrm
{u.c.p.}}{\longrightarrow}0$ which finishes the proof
of (\ref{sec3b}) for all $H$ with compact support.

\item[(iv)] Now, let $H\in C_p^1(\mathbb{R}^d)$ be arbitrary and $H_k$ be a
sequence of functions in $C_p^1(\mathbb{R}^d)$ with compact support
that converges pointwise to $H$ and fulfills $H=H_k$ on $[-k,k]^d$. In
view of step (i) it is enough to show that
\[
\lim_{k\to\infty} \limsup_{n \to\infty} \mathbb{E}
\biggl[\sup_{t\leq T} \biggl|\sqrt{n} \pmatrix{n
\cr
d}^{-1}
\sum_{\mathbf{i}\in
\mathcal{A}_t^n(d)} \nabla(H-H_k) \bigl(
\alpha_{\mathbf{i}}^n\bigr) \bigl(\sqrt{n}\Delta_{\mathbf{i}}^n
X- \alpha_{\mathbf{i}}^n\bigr) \biggr| \biggr]=0.
\]
Since $H-H_k$ is of polynomial growth and by (\ref{est1}), we get
\begin{eqnarray*}
&&\mathbb{E} \biggl[\sup_{t\leq T}  \biggl|\sqrt{n} \pmatrix{n
\cr
d}^{-1} \sum_{\mathbf{i}\in\mathcal{A}_t^n(d)} \nabla
(H-H_k) \bigl(\alpha_{\mathbf{i}}^n\bigr) \bigl(
\sqrt{n}\Delta_{\mathbf{i}}^n X- \alpha_{\mathbf{i}}^n
\bigr) \biggr| \biggr]
\\
&&\qquad \leq C\sqrt{n} \pmatrix{n
\cr
d}^{-1} \sum
_{\mathbf{i}\in\mathcal
{A}_T^n(d)} \mathbb{E}\bigl[ \bigl\|\nabla(H-H_k) \bigl(
\alpha_{\mathbf{i}}^n\bigr) \bigr\| \bigl\| \sqrt{n}\Delta_{\mathbf{i}}^n
X- \alpha_{\mathbf{i}}^n \bigr\| \bigr]
\\
&&\qquad \leq C \pmatrix{n
\cr
d}^{-1} \sum_{\mathbf{i}\in\mathcal
{A}_T^n(d)}
\mathbb{E} \Biggl[ \Biggl(\sum_{l=1}^d
\mathbh{1}_{ \{
|\alpha_{i_l}^n|>k \}} \Biggr)^2 \bigl\| \nabla(H-H_k)
\bigl(\alpha_{\mathbf{i}}^n\bigr) \bigr\|^2
\Biggr]^{1/2}\leq\frac{C}{k},
\end{eqnarray*}
which finishes the proof.
\end{longlist}

\subsection{Proof of (\texorpdfstring{\protect\ref{riemannphi}}{28})}\label{subsecsub33}
We can write
\[
U(H)_t = \int_{[0,t]^d} \int_{\mathbb{R}^d}
H(\mathbf{x}) \varphi_{\sigma_{s_1}}(x_1) \cdots
\varphi_{\sigma_{s_d}}(x_d) \,d\mathbf{x} \,d\mathbf{s}.
\]
We also have
\[
\widebar{F}_n'(t,x)=\int_0^{{[nt]}/{n}}
\varphi_{\sigma
_{{[ns]}/{n}}}(x) \,ds,
\]
where $\widebar{F}_n'(t,x)$ denotes the Lebesgue density in $x$ of
$\widebar{F}_n(t,x)$ defined at (\ref{riemannphi}). So we need to
show that $P^n(H) \stackrel{\mathrm{u.c.p.}}{\longrightarrow}0$, where
\begin{eqnarray*}
P_t^n(H)&:=& \sqrt{n}\int_{[0,t]^d} \int
_{\mathbb{R}^d} H(\mathbf{x})
\bigl( \varphi_{\sigma_{s_1}}(x_1) \cdots
\varphi_{\sigma_{s_d}}(x_d)
\\
&&\hspace*{87pt}{} -\varphi_{\sigma_{{[ns_1]}/{n}}}(x_1)
\cdots\varphi_{\sigma_{{[ns_d]}/{n}}}(x_d) \bigr) \,d\mathbf{x}
\,d\mathbf{s}.
\end{eqnarray*}
As previously we show the result first for $H$ with compact support.
\begin{longlist}[(ii)]
\item[(i)] Let the support of $H$ be contained in $[-k,k]^d$.
From \cite{BGJPS}, Section~8, we know that, for fixed $x \in\mathbb
{R}$, it holds that
%
%
\begin{equation}
\label{fcon2} \sqrt{n}\int_0^{t} \bigl(
\varphi_{\sigma_{s}}(x)-\varphi_{\sigma
_{{[ns]}/{n}}}(x) \bigr)\,ds \stackrel{\mathrm
{u.c.p.}} {\longrightarrow}0.
\end{equation}
Also, with $\rho(z,x):=\varphi_z(x)$ we obtain, for $x,y \in[-k,k]$,
\begin{eqnarray*}
&&\biggl|\int_0^t \bigl(\varphi_{\sigma_{s}}(x)-
\varphi_{\sigma_{{[ns]}/{n}}}(x)\bigr)-\bigl(\varphi_{\sigma_{s}}(y)-
\varphi_{\sigma_{{[ns]}/{n}}}(y)\bigr) \,ds \biggr|
\\
&&\qquad \leq\int_0^t \bigl| \partial_1 \rho(
\xi_s,x) (\sigma_{s}-\sigma_{{[ns]}/{n}})-
\partial_1 \rho\bigl(\xi_s',y\bigr) (
\sigma_{s}-\sigma_{{[ns]}/{n}}) \bigr| \,ds
\\
&&\qquad \leq\int_0^t \bigl|\partial_{11} \rho
\bigl(\xi_s'',\eta_s\bigr)
\bigl(\xi_s-\xi_s'\bigr)+
\partial_{21}\rho\bigl(\xi_s'',
\eta_s\bigr) (x-y) \bigr| |\sigma_{s}-\sigma_{{[ns]}/{n}} |
\,ds
\\
&&\qquad \leq C\int_0^t |\sigma_{s}-
\sigma_{{[ns]}/{n}}|^2+|\sigma_{s}-
\sigma_{{[ns]}/{n}}| |y-x| \,ds,
\end{eqnarray*}
where $\xi_s, \xi_s', \xi_s''$ are between $\sigma_{s}$ and $\sigma
_{[ns]/n}$ and $\eta_s$ is between $x$ and $y$. Now, let $Z_m= \{
jk/m | j=-m,\ldots,m \}$. Then, we get
\begin{eqnarray*}
\sup_{t \leq T} \bigl|P_t^n(H)\bigr| &\leq&
C_T\sup_{t \leq T} \sqrt{n} \int_{[-k,k]}
\biggl| \int_0^t \varphi_{\sigma_{s}}(x)-
\varphi_{\sigma
_{{[ns]}/{n}}}(x) \,ds \biggr|\,dx
\\
&\leq& C_T \sup_{t \leq T} \sup_{x \in[-k,k]}
\sqrt{n} \biggl| \int_0^t \varphi_{\sigma_{s}}(x)-
\varphi_{\sigma_{{[ns]}/{n}}}(x) \,ds \biggr|
\\
&\leq& C_T \sup_{t \leq T} \max_{x \in Z_m}
\sqrt{n} \biggl| \int_0^t \varphi_{\sigma_{s}}(x)-
\varphi_{\sigma_{{[ns]}/{n}}}(x) \,ds \biggr|
\\
&&{}+C_T\sqrt{n}\int_0^T
\biggl(|\sigma_{s}-\sigma_{{[ns]}/{n}}|^2+
\frac{k}{m}|\sigma_{s}-\sigma_{{[ns]}/{n}}|\biggr) \,ds
\\
&\leq& C_T \sum_{x \in Z_m} \sup
_{t \leq T} \sqrt{n} \biggl| \int_0^t
\varphi_{\sigma_{s}}(x)-\varphi_{\sigma_{{[ns]}/{n}}}(x) \,ds \biggr|
\\
&&{}+C_T\sqrt{n}\int_0^T
\biggl(|\sigma_{s}-\sigma_{{[ns]}/{n}}|^2+
\frac{k}{m}|\sigma_{s}-\sigma_{{[ns]}/{n}}|\biggr) \,ds.
\end{eqnarray*}
Observe that, for fixed $m$, the first summand converges in probability
to $0$ by (\ref{fcon2}). By the It\^o isometry and (\ref{probound2})
we get for the expectation of the second summand,
\begin{eqnarray*}
&&\mathbb{E} \biggl[ \sqrt{n}\int_0^T \biggl(|
\sigma_{s}-\sigma_{{[ns]}/{n}}|^2+\frac{k}{m}|
\sigma_{s}-\sigma_{{[ns]}/{n}}| \biggr) \,ds \biggr]
\\
&&\qquad =\sqrt{n} \int_0^T \mathbb{E} \biggl[|
\sigma_{s}-\sigma_{{[ns]}/{n}}|^2+\frac{k}{m}|
\sigma_{s}-\sigma_{{[ns]}/{n}}| \biggr] \,ds \le C_T
\biggl( \frac{1}{\sqrt{n}}+\frac{1}{m} \biggr).
\end{eqnarray*}
Thus, by choosing $m$ large enough and then letting $n$ go to infinity,
we get $P_t^n(H) \stackrel{\mathrm{u.c.p.}}{\longrightarrow}0$.

\item[(ii)] Now let $H \in C_p^1(\mathbb{R}^d)$ and $H_k$ be an approximating
sequence of functions in $C_p^1(\mathbb{R}^d)$ with compact support
and $H=H_k$ on $[-k,k]^d$. Observe that, for \mbox{$\mathbf{x},\mathbf
{s}\in\mathbb{R}^d$}, we obtain by the mean value theorem that
\begin{eqnarray*}
&& \mathbb{E} \bigl[ \bigl| \varphi_{\sigma_{s_1}}(x_1) \cdots\varphi
_{\sigma_{s_d}}(x_d)-\varphi_{\sigma_{{[ns_1]}/{n}}}(x_1)
\cdots\varphi_{\sigma_{{[ns_d]}/{n}}}(x_d) \bigr| \bigr]
\\
&&\qquad \leq\psi(\mathbf{x}) \sum_{i=1}^d
\mathbb{E}|\sigma_{s_i}-\sigma_{{[ns_i]}/{n}}| \leq\frac{C}{\sqrt{n}}
\psi(\mathbf{x}),
\end{eqnarray*}
where the function $\psi$ is exponentially decaying at $\pm\infty$. Thus
\begin{eqnarray*}
&& \lim_{k \to\infty} \limsup_{n \to\infty} \mathbb{E}
\Bigl[\sup_{t\leq T} \bigl|P_t^n
(H)-P_t^n(H_k)\bigr| \Bigr]
\\
&&\qquad \leq C_T \lim_{k \to\infty} \limsup
_{n \to\infty} \int_{\mathbb
{R}^d} \bigl|(H-H_k) (
\mathbf{x})\bigr| \psi(\mathbf{x}) \,d\mathbf{x} =0,
\end{eqnarray*}
which finishes the proof of (\ref{riemannphi}).
\end{longlist}

\section*{Acknowledgment}
We would like to thank Herold Dehling for his helpful comments.



%

\printaddresses


\begin{thebibliography}{27}
\bibitem{AE}
%
\begin{barticle}[mr]
\bauthor{\bsnm{Aldous},~\bfnm{D.~J.}\binits{D.~J.}} \AND
\bauthor{\bsnm{Eagleson},~\bfnm{G.~K.}\binits{G.~K.}}
(\byear{1978}).
\btitle{On mixing and stability of limit theorems}.
\bjournal{Ann. Probab.}
\bvolume{6}
\bpages{325--331}.
\bid{mr={0517416}}
\end{barticle}
%
\bptok{imsref}%
\endbibitem

\bibitem{BGJPS}
%
\begin{bincollection}[mr]
\bauthor{\bsnm{Barndorff-Nielsen},~\bfnm{Ole~E.}\binits{O.~E.}},
\bauthor{\bsnm{Graversen},~\bfnm{Svend~Erik}\binits{S.~E.}},
\bauthor{\bsnm{Jacod},~\bfnm{Jean}\binits{J.}},
\bauthor{\bsnm{Podolskij},~\bfnm{Mark}\binits{M.}} \AND
\bauthor{\bsnm{Shephard},~\bfnm{Neil}\binits{N.}}
(\byear{2006}).
\btitle{A central limit theorem for realised power and bipower
variations of continuous semimartingales}.
In \bbooktitle{From Stochastic Calculus to Mathematical Finance.
Festschrift in Honour of A.~N.~Shiryaev}
(\beditor{\bfnm{Yu.}\binits{Y.}~\bsnm{Kabanov}},
\beditor{\bfnm{R.}\binits{R.}~\bsnm{Liptser}} \AND
\beditor{\bfnm{J.}\binits{J.}~\bsnm{Stoyanov}}, eds.)
\bpages{33--68}.
\bpublisher{Springer},
\blocation{Heidelberg}.
\bid{doi={10.1007/978-3-540-30788-4_3}, mr={2233534}}
\end{bincollection}
%
\bptok{imsref}%
\endbibitem

\bibitem{BZ}
%
\begin{barticle}[mr]
\bauthor{\bsnm{Beutner},~\bfnm{Eric}\binits{E.}} \AND
\bauthor{\bsnm{Z{\"a}hle},~\bfnm{Henryk}\binits{H.}}
(\byear{2012}).
\btitle{Deriving the asymptotic distribution of $U$- and \mbox{$V$-}statistics
of dependent data using weighted empirical processes}.
\bjournal{Bernoulli}
\bvolume{18}
\bpages{803--822}.
\bid{doi={10.3150/11-BEJ358}, issn={1350-7265}, mr={2948902}}
\end{barticle}
%
\bptok{imsref}%
\endbibitem

\bibitem{BBD1}
%
\begin{barticle}[mr]
\bauthor{\bsnm{Borovkova},~\bfnm{S.}\binits{S.}},
\bauthor{\bsnm{Burton},~\bfnm{R.}\binits{R.}} \AND
\bauthor{\bsnm{Dehling},~\bfnm{H.}\binits{H.}}
(\byear{1999}).
\btitle{Consistency of the {T}akens estimator for the correlation dimension}.
\bjournal{Ann. Appl. Probab.}
\bvolume{9}
\bpages{376--390}.
\bid{doi={10.1214/aoap/1029962747}, issn={1050-5164}, mr={1687339}}
\end{barticle}
%
\bptok{imsref}%
\endbibitem

\bibitem{BBD2}
%
\begin{barticle}[mr]
\bauthor{\bsnm{Borovkova},~\bfnm{Svetlana}\binits{S.}},
\bauthor{\bsnm{Burton},~\bfnm{Robert}\binits{R.}} \AND
\bauthor{\bsnm{Dehling},~\bfnm{Herold}\binits{H.}}
(\byear{2001}).
\btitle{Limit theorems for functionals of mixing processes with
applications to $U$-statistics and dimension estimation}.
\bjournal{Trans. Amer. Math. Soc.}
\bvolume{353}
\bpages{4261--4318}.
\bid{doi={10.1090/S0002-9947-01-02819-7}, issn={0002-9947}, mr={1851171}}
\end{barticle}
%
\bptok{imsref}%
\endbibitem

\bibitem{CH}
%
\begin{bbook}[author]
\bauthor{\bsnm{Cs{\"o}rg{\'o}},~\bfnm{M.}\binits{M.}} \AND
\bauthor{\bsnm{Horv{\'a}th},~\bfnm{L.}\binits{L.}}
(\byear{1997}).
\btitle{Limit Theorems in Change-point Analysis}.
\bpublisher{Wiley},
\blocation{New York}.
\end{bbook}
%
\bptok{imsref}%
\endbibitem

\bibitem{David1968}
%
\begin{barticle}[author]
\bauthor{\bsnm{David},~\bfnm{H.~A.}\binits{H.~A.}}
(\byear{1968}).
\btitle{Gini's mean difference rediscovered}.
\bjournal{Biometrika}
\bvolume{55}
\bpages{573--575}.
\end{barticle}
%
\bptok{imsref}%
\endbibitem

\bibitem{DRT}
%
\begin{barticle}[mr]
\bauthor{\bsnm{Dehling},~\bfnm{Herold}\binits{H.}},
\bauthor{\bsnm{Rooch},~\bfnm{Aeneas}\binits{A.}} \AND
\bauthor{\bsnm{Taqqu},~\bfnm{Murad~S.}\binits{M.~S.}}
(\byear{2013}).
\btitle{Nonparametric change-point tests for long-range dependent data}.
\bjournal{Scand. J. Stat.}
\bvolume{40}
\bpages{153--173}.
\bid{doi={10.1111/j.1467-9469.2012.00799.x}, issn={0303-6898}, mr={3024037}}
\end{barticle}
%
\bptok{imsref}%
\endbibitem

\bibitem{DT}
%
\begin{barticle}[mr]
\bauthor{\bsnm{Dehling},~\bfnm{Herold}\binits{H.}} \AND
\bauthor{\bsnm{Taqqu},~\bfnm{Murad~S.}\binits{M.~S.}}
(\byear{1989}).
\btitle{The empirical process of some long-range dependent sequences
with an application to $U$-statistics}.
\bjournal{Ann. Statist.}
\bvolume{17}
\bpages{1767--1783}.
\bid{doi={10.1214/aos/1176347394}, issn={0090-5364}, mr={1026312}}
\end{barticle}
%
\bptok{imsref}%
\endbibitem

\bibitem{DTb}
%
\begin{barticle}[mr]
\bauthor{\bsnm{Dehling},~\bfnm{Herold}\binits{H.}} \AND
\bauthor{\bsnm{Taqqu},~\bfnm{Murad~S.}\binits{M.~S.}}
(\byear{1991}).
\btitle{Bivariate symmetric statistics of long-range dependent observations}.
\bjournal{J. Statist. Plann. Inference}
\bvolume{28}
\bpages{153--165}.
\bid{doi={10.1016/0378-3758(91)90031-9}, issn={0378-3758}, mr={1115815}}
\end{barticle}
%
\bptok{imsref}\vadjust{\goodbreak}%
\endbibitem

\bibitem{DK}
%
\begin{barticle}[mr]
\bauthor{\bsnm{Denker},~\bfnm{Manfred}\binits{M.}} \AND
\bauthor{\bsnm{Keller},~\bfnm{Gerhard}\binits{G.}}
(\byear{1983}).
\btitle{On $U$-statistics and v. {M}ises' statistics for weakly dependent
processes}.
\bjournal{Z. Wahrsch. Verw. Gebiete}
\bvolume{64}
\bpages{505--522}.
\bid{doi={10.1007/BF00534953}, issn={0044-3719}, mr={0717756}}
\end{barticle}
%
\bptok{imsref}%
\endbibitem

\bibitem{DP}
%
\begin{barticle}[mr]
\bauthor{\bsnm{Dette},~\bfnm{Holger}\binits{H.}} \AND
\bauthor{\bsnm{Podolskij},~\bfnm{Mark}\binits{M.}}
(\byear{2008}).
\btitle{Testing the parametric form of the volatility in continuous
time diffusion models---a stochastic process approach}.
\bjournal{J. Econometrics}
\bvolume{143}
\bpages{56--73}.
\bid{doi={10.1016/j.jeconom.2007.08.002}, issn={0304-4076}, mr={2384433}}
\end{barticle}
%
\bptok{imsref}%
\endbibitem

\bibitem{DPV}
%
\begin{barticle}[mr]
\bauthor{\bsnm{Dette},~\bfnm{Holger}\binits{H.}},
\bauthor{\bsnm{Podolskij},~\bfnm{Mark}\binits{M.}} \AND
\bauthor{\bsnm{Vetter},~\bfnm{Mathias}\binits{M.}}
(\byear{2006}).
\btitle{Estimation of integrated volatility in continuous-time
financial models with applications to goodness-of-fit testing}.
\bjournal{Scand. J. Stat.}
\bvolume{33}
\bpages{259--278}.
\bid{doi={10.1111/j.1467-9469.2006.00479.x}, issn={0303-6898}, mr={2279642}}
\end{barticle}
%
\bptok{imsref}%
\endbibitem

\bibitem{GT}
%
\begin{barticle}[mr]
\bauthor{\bsnm{Giraitis},~\bfnm{Liudas}\binits{L.}} \AND
\bauthor{\bsnm{Taqqu},~\bfnm{Murad~S.}\binits{M.~S.}}
(\byear{1997}).
\btitle{Limit theorems for bivariate {A}ppell polynomials. {I}.~{C}entral limit theorems}.
\bjournal{Probab. Theory Related Fields}
\bvolume{107}
\bpages{359--381}.
\bid{doi={10.1007/s004400050089}, issn={0178-8051}, mr={1440137}}
\end{barticle}
%
\bptok{imsref}%
\endbibitem

\bibitem{Hoe}
%
\begin{barticle}[mr]
\bauthor{\bsnm{Hoeffding},~\bfnm{Wassily}\binits{W.}}
(\byear{1948}).
\btitle{A class of statistics with asymptotically normal distribution}.
\bjournal{Ann. Math. Stat.}
\bvolume{19}
\bpages{293--325}.
\bid{issn={0003-4851}, mr={0026294}}
\end{barticle}
%
\bptok{imsref}%
\endbibitem

\bibitem{HW}
%
\begin{barticle}[mr]
\bauthor{\bsnm{Hsing},~\bfnm{Tailen}\binits{T.}} \AND
\bauthor{\bsnm{Wu},~\bfnm{Wei~Biao}\binits{W.~B.}}
(\byear{2004}).
\btitle{On weighted $U$-statistics for stationary processes}.
\bjournal{Ann. Probab.}
\bvolume{32}
\bpages{1600--1631}.
\bid{doi={10.1214/009117904000000333}, issn={0091-1798}, mr={2060311}}
\end{barticle}
%
\bptok{imsref}%
\endbibitem

\bibitem{J3}
%
\begin{bincollection}[mr]
\bauthor{\bsnm{Jacod},~\bfnm{Jean}\binits{J.}}
(\byear{1997}).
\btitle{On continuous conditional {G}aussian martingales and stable
convergence in law}.
In \bbooktitle{S\'eminaire de {P}robabilit\'es, {XXXI}}.
\bseries{Lecture Notes in Math.}
\bvolume{1655}
\bpages{232--246}.
\bpublisher{Springer},
\blocation{Berlin}.
\bid{doi={10.1007/BFb0119308}, mr={1478732}}
\end{bincollection}
%
\bptok{imsref}%
\endbibitem

\bibitem{J2}
%
\begin{barticle}[mr]
\bauthor{\bsnm{Jacod},~\bfnm{Jean}\binits{J.}}
(\byear{2008}).
\btitle{Asymptotic properties of realized power variations and related
functionals of semimartingales}.
\bjournal{Stochastic Process. Appl.}
\bvolume{118}
\bpages{517--559}.
\bid{doi={10.1016/j.spa.2007.05.005}, issn={0304-4149}, mr={2394762}}
\end{barticle}
%
\bptok{imsref}%
\endbibitem

\bibitem{JPV}
%
\begin{barticle}[mr]
\bauthor{\bsnm{Jacod},~\bfnm{Jean}\binits{J.}},
\bauthor{\bsnm{Podolskij},~\bfnm{Mark}\binits{M.}} \AND
\bauthor{\bsnm{Vetter},~\bfnm{Mathias}\binits{M.}}
(\byear{2010}).
\btitle{Limit theorems for moving averages of discretized processes
plus noise}.
\bjournal{Ann. Statist.}
\bvolume{38}
\bpages{1478--1545}.
\bid{doi={10.1214/09-AOS756}, issn={0090-5364}, mr={2662350}}
\end{barticle}
%
\bptok{imsref}%
\endbibitem

\bibitem{JP}
%
\begin{bbook}[mr]
\bauthor{\bsnm{Jacod},~\bfnm{Jean}\binits{J.}} \AND
\bauthor{\bsnm{Protter},~\bfnm{Philip}\binits{P.}}
(\byear{2012}).
\btitle{Discretization of Processes}.
\bpublisher{Springer},
\blocation{Heidelberg}.
\bid{doi={10.1007/978-3-642-24127-7}, mr={2859096}}
\end{bbook}
%
\bptok{imsref}%
\endbibitem

\bibitem{JS}
%
\begin{bbook}[mr]
\bauthor{\bsnm{Jacod},~\bfnm{Jean}\binits{J.}} \AND
\bauthor{\bsnm{Shiryaev},~\bfnm{Albert~N.}\binits{A.~N.}}
(\byear{2003}).
\btitle{Limit Theorems for Stochastic Processes},
\bedition{2nd} ed.
\bpublisher{Springer},
\blocation{Berlin}.
\bid{mr={1943877}}
\end{bbook}
%
\bptok{imsref}%
\endbibitem

\bibitem{KP}
%
\begin{barticle}[mr]
\bauthor{\bsnm{Kinnebrock},~\bfnm{Silja}\binits{S.}} \AND
\bauthor{\bsnm{Podolskij},~\bfnm{Mark}\binits{M.}}
(\byear{2008}).
\btitle{A note on the central limit theorem for bipower variation of
general functions}.
\bjournal{Stochastic Process. Appl.}
\bvolume{118}
\bpages{1056--1070}.
\bid{doi={10.1016/j.spa.2007.07.009}, issn={0304-4149}, mr={2418258}}
\end{barticle}
%
\bptok{imsref}%
\endbibitem

\bibitem{KoroljukBorovskich1994}
%
\begin{bbook}[mr]
\bauthor{\bsnm{Koroljuk},~\bfnm{V.~S.}\binits{V.~S.}} \AND
\bauthor{\bsnm{Borovskich},~\bfnm{Yu.~V.}\binits{Yu.~V.}}
(\byear{1994}).
\btitle{Theory of $U$-Statistics}.
\bpublisher{Kluwer},
\blocation{Dordrecht}.
\bid{mr={1472486}}
\end{bbook}
%
\bptok{imsref}%
\endbibitem

\bibitem{LL}
%
\begin{barticle}[mr]
\bauthor{\bsnm{L{\'e}vy-Leduc},~\bfnm{C.}\binits{C.}},
\bauthor{\bsnm{Boistard},~\bfnm{H.}\binits{H.}},
\bauthor{\bsnm{Moulines},~\bfnm{E.}\binits{E.}},
\bauthor{\bsnm{Taqqu},~\bfnm{M.~S.}\binits{M.~S.}} \AND
\bauthor{\bsnm{Reisen},~\bfnm{V.~A.}\binits{V.~A.}}
(\byear{2011}).
\btitle{Asymptotic properties of $U$-processes under long-range dependence}.
\bjournal{Ann. Statist.}
\bvolume{39}
\bpages{1399--1426}.
\bid{doi={10.1214/10-AOS867}, issn={0090-5364}, mr={2850207}}
\end{barticle}
%
\bptok{imsref}%
\endbibitem

\bibitem{PV}
%
\begin{barticle}[mr]
\bauthor{\bsnm{Podolskij},~\bfnm{Mark}\binits{M.}} \AND
\bauthor{\bsnm{Vetter},~\bfnm{Mathias}\binits{M.}}
(\byear{2010}).
\btitle{Understanding limit theorems for semimartingales: A short survey}.
\bjournal{Stat. Neerl.}
\bvolume{64}
\bpages{329--351}.
\bid{doi={10.1111/j.1467-9574.2010.00460.x}, issn={0039-0402}, mr={2683464}}
\end{barticle}
%
\bptok{imsref}%
\endbibitem

\bibitem{R}
%
\begin{barticle}[mr]
\bauthor{\bsnm{R{\'e}nyi},~\bfnm{Alfr{\'e}d}\binits{A.}}
(\byear{1963}).
\btitle{On stable sequences of events}.
\bjournal{Sankhy\=a Ser. A}
\bvolume{25}
\bpages{293--302}.
\bid{issn={0581-572X}, mr={0170385}}
\end{barticle}
%
\bptok{imsref}%
\endbibitem

\bibitem{Yitzhaki1982}
%
\begin{barticle}[author]
\bauthor{\bsnm{Yitzhaki},~\bfnm{S.}\binits{S.}}
(\byear{1982}).
\btitle{Stochastic dominance, mean variance, and {G}ini's mean difference}.
\bjournal{Am. Econ. Rev.}
\bvolume{71}
\bpages{178--185}.
\end{barticle}
%
\bptok{imsref}%
\endbibitem

\end{thebibliography}
\end{document}